# THE PALM MEASURE AND THE VORONOI TESSELLATION FOR THE GINIBRE PROCESS

By André Goldman

*Université Claude Bernard Lyon 1*

We prove that the Palm measure of the Ginibre process is obtained by removing a Gaussian distributed point from the process and adding the origin. We obtain also precise formulas describing the law of the typical cell of Ginibre–Voronoi tessellation. We show that near the germs of the cells a more important part of the area is captured in the Ginibre–Voronoi tessellation than in the Poisson–Voronoi tessellation. Moment areas of corresponding subdomains of the cells are explicitly evaluated.

**1. Introduction and statement of the main results.** The Poisson–Voronoi tessellation is a very popular model of stochastic geometry. This is mainly due to its large range of applicability: crystallography [7], astrophysics [32] and telecommunications [1], to mention only a few. This is also due to the simplicity of the simulation procedures [11, 13, 31], and to the fact that several theoretical results related to its geometrical characteristics are available [4, 9, 10, 21, 22]. An extensive list of the areas in which these tessellations have been used can be found in [23, 29]. Nevertheless, the other side of the picture is that the comparative triviality of this model makes it inappropriate to describe precisely some natural phenomena. Hence, it seems both interesting and useful to explore other random point processes and their Voronoi tessellations. For instance, Le Caer and Ho [16] describe, by means of Monte Carlo simulations, statistical properties of the Voronoi tessellation associated to the Ginibre process of eigenvalues of random complex Gaussian matrices [8] (see also [26]). The idea behind their study is that the repulsive character of the distribution of random points makes the cells more regular. Consequently, the associated tessellation fits better than the Poisson–Voronoi one, as in, for example, the structure of the cells of biological tissues. We recall that the Ginibre process [20, 28], is a determinantal process $\phi \subset R^2$, both isotropic and ergodic with respect to the translations of the plane $R^2 = \mathbb{C}$, with the integral kernel,

$$K(z_1, z_2) = (1/\pi) e^{z_1 \overline{z}_2} \exp(-(1/2)(|z_1|^2 + |z_2|^2)), \qquad (z_1, z_2) \in \mathbb{C}^2. \tag{1}$$

It is also pertinent to consider the full class of determinantal processes $\phi^{\star\alpha}$ related to the kernels, $K_{\star\alpha}(z_1, z_2) = (1/\pi) e^{(1/\alpha) z_1 \overline{z}_2} \exp(-(1/2\alpha)(|z_1|^2 + |z_2|^2))$, $(z_1, z_2) \in \mathbb{C}^2$, with $0 < \alpha < 1$. The process $\phi^{\star\alpha}$ can be obtained by deleting, independently and with probability $1 - \alpha$, each point of the Ginibre process $\phi$ and then applying the homothety of ratio $\sqrt{\alpha}$ to the remaining points in order to restore the intensity of the process $\phi$. Besides, it is easy to verify that $\phi^{\star\alpha}$ converges in law when $\alpha \to 0$ to the Poisson process. In other words, the processes $\phi^{\star\alpha}$ constitute an intermediate class between a Poisson process and a Ginibre process. In order to challenge the classical Poisson–Voronoi model, it is necessary to have some theoretical knowledge about geometric characteristics of Ginibre–Voronoi tessellations. The main tool for this is the notion of a typical cell in the Palm sense [22]. To explain this notion, we introduce, for a general stationary process $\psi$, the notation,

$$\psi_0 = (\psi \mid 0 \in \psi) \setminus \{0\}.$$

The typical cell of $\psi$ is

$$\mathcal{C} = \{z \in \mathbb{C}; \forall u \in \psi_0, |z| \le |z - u|\}.$$

When $\psi$ is ergodic, the laws of the geometric characteristics of the typical cell coincide (see [3, 5, 6]) with the empirical distributions of the corresponding characteristics associated to the Voronoi tessellation, $\{C(u, \psi); u \in \psi\}$, whose cells are

$$C(u, \psi) = \{z \in \mathbb{C}; \forall v \in \psi, |z - u| \le |z - v|\}, \qquad u \in \psi.$$

If $\psi$ is a Poisson stationary process, then the Slivnyak formula [22] states that, for every finite set $S \subset \mathbb{C}$,

$$(\psi \mid S \subset \psi) \setminus S \stackrel{\text{law}}{=} \psi.$$

Hence, in this case, the Palm measure of $\psi$ is the law of the process $\psi \cup \{0\}$ obtained by adding the origin to $\psi$. For every determinantal process $\psi$, a result obtained by Shirai and Takahashi [27] states that $\psi_0$ is determinantal as well. It follows that in the Ginibre case the process $\phi_0$ is determinantal with the kernel,

$$K_0(z_1, z_2) = (1/\pi)(e^{z_1 \overline{z}_2} - 1) \exp(-(1/2)(|z_1|^2 + |z_2|^2)), \qquad (z_1, z_2) \in \mathbb{C}^2, \tag{2}$$

and that

$$\mathcal{C} \stackrel{\text{law}}{=} C(0, \phi_0).$$

Note that the process $\phi_0$ is nonstationary.

Our first main result is that $\phi_0$ can be obtained from $\phi$, simply by deleting one point. A more precise statement follows.

THEOREM 1. *There exists a Gaussian-centered random variable $Z$, such that $E|Z|^2 = 1$ and*

$$\phi \overset{\text{law}}{=} \phi_0 \cup \{Z\}, \qquad \phi_0 \cap \{Z\} = \varnothing.$$

Theorem 1 tells us that there exists a version of the Ginibre process $\phi$ such that the Palm measure of $\phi$ is the law of the process obtained by removing from $\phi$ a Gaussian-distributed point and then adding the origin. As an intermediate step on our way to further results, consider a locally compact Hausdorff space $\mathbb{E}$ with a countable basis and a reference Radon measure $\lambda$, and a general stationary determinantal process $\psi \subset \mathbb{E}$ with kernel $K$ defined on $\mathbb{E}^2$. We introduce the following conditions.

CONDITION I. The measure $\lambda$ has full support, that is, for every open set $U \subset \mathbb{E}$, $\lambda(U)$ is positive.

CONDITION A. The kernel $K$ is a continuous function on $\mathbb{E}^2$.

CONDITION B. For every bounded Borel set $A \subset \mathbb{E}$, all the eigenvalues of the operator $K_A$, acting on $L^2(A, \lambda)$, lie in the interval $[0, 1[$.

Our intermediate result is the following.

THEOREM 2. *Assume that $\mathbb{E}$ and $\psi$ satisfy Conditions I, A and B. Then the process $\psi_0$ is stochastically dominated by the process $\psi$.*

More generally we prove that:

THEOREM 3. *Assume that $\mathbb{E}$ satisfies Condition I and consider two kernels $K$ and $L$ satisfying Conditions A and B above. Denote by $\psi$ the determinantal process associated to the kernel $K$ and by $\varphi$ the process associated to the kernel $L$. Suppose that $K \geq L$ in the Loewner order. Then the process $\psi$ dominates stochastically the process $\varphi$.*

Recall that $K \geq L$ in the Loewner order if $K - L$ is a positive semidefinite operator. For the kernels $K$ and $K_0$ defined by formulas (1) and (2), we have, obviously, $K \geq K_0$ thus Theorem 3 implies Theorem 2.

The proof of Theorem 3 rests on an explicit description of the marginal laws of the process $\psi$, obtained in Section 1, which allows us to use a similar result proved by Lyons [19] (for commuting operators) and for Borcea, Branden and Liggett [2] (without this restriction) in the discrete determinantal process setting (see also errata to [19] on Russel Lyon's website).

Thanks to the characterization of the Palm measure, we obtain, following [4], precise formulas (see Section 5) which describe the law of the typical cell of the Ginibre–Voronoi tessellation. The integrals involved are rather awkward; this should not be a surprise, since this is already the case for the Poisson–Voronoi typical cell [4].

In the last part of this work we compare the moments of the areas of Poisson–Voronoi and Ginibre–Voronoi cells. We show that near the germs of the cells, a more important part of the area is captured in the Ginibre case; farther from the germs of the cells, the situation is reversed. That is, roughly speaking, Ginibre cells are more stocky than Poisson cells.

To be more precise, we introduce some notation. Let $\mathcal{C}_p$ denote the typical cell of the Voronoi tessellation associated to a stationary Poisson process in $\mathbb{C}$ with the same intensity as the process $\phi$. For every positive $r$ and every $z \in \mathbb{C}$, let $B(z,r) \subset \mathbb{C}$ denote the disc centered at $z$ with radius $r$, and $B(r) = B(0,r)$. For every finite set $S \subset \mathbb{C}$,

$$D(S) = \bigcup_{z \in S} B(z, |z|).$$

Let $V(S)$ denote the area of $D(S)$. For every Borel set $A \subset \mathbb{C}$, and every positive integer $k$ introduce

$$V^k(A) = \left[\int_A dz\right]^k,$$

where, for $z = x + \mathrm{i}y$ in $\mathbb{C}$ with $(x,y) \in \mathbb{R}^2$, one sets $dz = dx\,dy$. Finally, for $z = (z_1, \ldots, z_k) \in \mathbb{C}^k$, one sets $dz = dz_1 \cdots dz_k$.

THEOREM 4. *Let $k$ denote a positive integer.*

(a) *When $r \to 0$,*

$$EV^k(\mathcal{C} \cap B(r)) = EV^k(\mathcal{C}_p \cap B(r))(1 + r^2 W_k + o(r^2)) \tag{3}$$

*with*

$$W_k = \frac{1}{\pi^{k+1}} \int_{B(1)^k} V(z)\,dz.$$

(b) *For every positive $R$,*

$$EV^k(\mathcal{C} \setminus B(R)) \le EV^k(\mathcal{C}_p \setminus B(R)) \cdot e^{(3/2) - J(R)} \tag{4}$$

*with*

$$J(R) = \frac{1}{2\pi^2} \int_{B(R)^2} e^{-|z_1 - z_2|^2}\,dz_1\,dz_2.$$

*Hence, there exists a positive constant $c$ such that, for every positive $R$,*

$$EV^k(\mathcal{C} \setminus B(R)) \le EV^k(\mathcal{C}_p \setminus B(R)) \cdot e^{(3/2)-cR^2}.$$

We are also interested in the location of the point $Z$ in Theorem 1 with respect to the process $\phi_0$. To describe this, let $\mathcal{N}(\phi)$ denote the set of points $z \in \phi$ such that the bisecting line of the segment $[0, z]$ intersects the boundary of the cell $C(0, \phi)$ where we recall that

$$C(0, \phi) = \{z \in \mathbb{C}; \forall u \in \phi, |z| \le |z - u|\}.$$

For every set $S \subset \mathbb{C}$, let

$$H(S) = (\pi R_{D(S)}(0,0) - 1) \prod_{n \ge 0} (1 - \alpha_n(S)), \tag{5}$$

where $R_{D(S)}$ is the resolvent kernel, and $\alpha_n(S)$ for $n \ge 0$, are the eigenvalues of the integral operator $K$ acting on the space $L^2(D(S), dz)$.

Theorem 5. *For every positive integer $k$,*

$$P\{Z \in \mathcal{N}(\phi)\} \ge \left[ \frac{\int_{\mathbb{C}^k} H(z)\, dz}{\int_{\mathbb{C}^k} \sqrt{H(z)}\, dz} \right]^2. \tag{6}$$

Taking $k = 1$, Theorem 5 yields

$$\begin{aligned} &P\{Z \in \mathcal{N}(\phi)\} \\ &\ge \left( \frac{1 - \int_0^{+\infty} \prod_{n \ge 0} (\Gamma(n+1, t)/n!)\, dt}{\int_0^{+\infty} \sqrt{(\sum_{n \ge 1} (t^n e^{-t} / \Gamma(n+1, t))) \prod_{n \ge 0} (\Gamma(n+1, t)/n!)}\, dt} \right)^2, \end{aligned} \tag{7}$$

where $\Gamma(n, t)$ denotes the incomplete gamma function, defined as

$$\Gamma(n, t) = \int_t^{+\infty} e^{-u} u^{n-1}\, du.$$

This implies the simpler bound,

$$P\{Z \in \mathcal{N}(\phi)\} \ge \tfrac{1}{16}.$$

Section 2 contains the necessary background; the key results are Propositions 3 and 12. In Section 3 we prove Theorem 3. Theorem 1 is proved in Section 4 as a consequence of Theorem 2, Strassen's classical result and because of the fact that the radial processes $|\phi|$ and $|\phi_0|$ are explicitly known. Unfortunately, the correlation between the process $\phi_0$ and the random point $Z$ is still unknown. Nevertheless, Theorem 5 gives some partial insight in this direction. Finally, we mention that we state our results for the Ginibre process, but that it is easy to deduce the corresponding formulations for the processes $\phi^{\star\alpha}$ with $0 < \alpha < 1$.

**2. Preliminaries.** Let $\psi$ denote a point process [28] on a locally compact Hausdorff space $\mathbb{E}$ with a countable basis. For every integer $k \geq 1$, let $\psi^{(k)} = \{\tilde{z} \subset \psi; |\tilde{z}| = k\}$, $\psi^{(1)} = \psi$, be the associated $k$-dimensional process. Let $\mu_k$ denote the corresponding intensity measure. This measure is defined as follows. Fix a set $\tilde{z} \in \psi^{(k)}$ and consider an arbitrary order $\tilde{z} = \{z_1, \dots, z_k\}$. For every permutation $\sigma \in \wp_k$ of the index set $\{1, \dots, k\}$, denote

$$z^\sigma = (z_{\sigma(1)}, \dots, z_{\sigma(k)}) \in \mathbb{E}^k.$$

For every Borel set $A \in \mathcal{B}(\mathbb{E}^k)$ of the space $\mathbb{E}^k$, the sum $\sum_{\sigma \in \wp_k} 1_A(z^\sigma)$ does not depend on the particular ordering of the set $\tilde{z}$. By summing for $\tilde{z} \in \psi^k$ and taking the expectation, we obtain

$$\mu_k(A) = \frac{1}{k!} E \sum_{\tilde{z} \in \psi^{(k)}} \sum_{\sigma \in \wp_k} 1_A(z^\sigma).$$

Consider the space $\mathcal{M}_\sigma(\mathbb{E})$ of the counting measures $\xi$ on $\mathbb{E}$ such that $\xi(A)$ is finite for all bounded (relatively compact) Borel sets $A \subset \mathbb{E}$, and let $\mathcal{F}$ be the smallest $\sigma$-algebra on $\mathcal{M}_\sigma(\mathbb{E})$ for which the map $\xi \mapsto \xi(A)$ is measurable for every bounded Borel set $A \subset \mathbb{E}$.

The point process $\psi$ can be thought of as the random measure

$$\xi = \sum_{z \in \psi} \delta_z$$

with values in the measurable space $(\mathcal{M}_\sigma(\mathbb{E}), \mathcal{F})$. Note that the space $\mathcal{M}_\sigma(\mathbb{E})$ endowed with vague topology is a Polish space and that the associated Borel $\sigma$-algebra coincides with the $\sigma$-algebra $\mathcal{F}$ (see [17] and [18]).

For every $k \geq 1$, the Campbell measure $C_k$ on $\mathbb{E}^k \times \mathcal{M}_\sigma(\mathbb{E})$ is

$$C_k(M) = \frac{1}{k!} E \sum_{\tilde{z} \in \psi^{(k)}} \sum_{\sigma \in \wp_k} 1_M(z^\sigma, \psi), \qquad M \in \mathcal{B}(\mathbb{E}^k) \otimes \mathcal{F},$$

where, as above, the sum $\sum_{\sigma \in \wp_k} 1_M(z^\sigma, \psi)$, $z^\sigma = (z_{\sigma(1)}, \dots, z_{\sigma(k)})$, does not depend on the particular ordering $\tilde{z} = \{z_1, \dots, z_k\}$ of the set $\tilde{z}$.

The disintegration of $C_k$ with respect to the measure $\mu_k$ gives, for $\mu_k$ almost every $z = (z_1, \dots, z_k) \in \mathbb{E}^k$, the law of the conditioned process $(\psi \mid \tilde{z} = \{z_1, \dots, z_k\} \in \psi^{(k)})$ (see [14]). Campbell formula reads as follows (see [14]). Assume that $f$ is a measurable positive function defined on $\mathbb{E}^k \times \mathcal{M}_\sigma(\mathbb{E})$, such that $f(z, \cdot) = f(z^\sigma, \cdot)$ for every permutation $\sigma \in \wp_k$ and for $\mu_k$ almost every $z \in \mathbb{E}^k$, thus $f$ defines a function acting on sets $\tilde{z} \in \psi^{(k)}$ by $f(\tilde{z}, \cdot) = f((z_1, \dots, z_k), \cdot)$ where $\tilde{z} = \{z_1, \dots, z_k\}$ is an arbitrary ordering. Then,

$$\begin{aligned} & E \sum_{\tilde{z} \in \psi^{(k)}} f(\tilde{z}, \psi) \\ & \quad = \int E f((z_1, \dots, z_k), (\psi \mid \{z_1, \dots, z_k\} \subset \psi)) \, d\mu_k(z_1, \dots, z_k). \end{aligned} \tag{8}$$

Let $\lambda$ denote a Radon measure on $\mathbb{E}$, that is, a Borel measure such that $\lambda(A)$ is finite for every compact set $A \subset \mathbb{E}$. In this paper $\lambda$ will be the Lebesgue measure on $\mathbb{C} = \mathbb{R}^2$ or the standard counting measure on a finite discrete set. The point process $\psi$ is determinantal if the following properties hold.

1. For every $k \geq 1$, $\mu_k$ is absolutely continuous with respect to the product measure $\lambda^k$ on $\mathbb{E}^k$, that is, there exists a density $\rho_k$ such that
$$d\mu_k = \rho_k \, d\lambda^k.$$
The density $\rho_k$ is called the correlation function.
2. There exists a kernel $K : \mathbb{E} \times \mathbb{E} \to \mathbb{C}$ which defines a self-adjoint, locally trace-class operator, such that, for every $z = (z_1, \dots, z_k)$ in $\mathbb{E}^k$,

$$\rho_k(z) = \frac{1}{k!} \det(K(z_i, z_j))_{1 \leq i,j \leq k}. \tag{9}$$

We use Fredholm notation; hence for every $k \geq 1$ and every $u = (u_1, \dots, u_k)$ and $v = (v_1, \dots, v_k)$ in $\mathbb{E}^k$,
$$K\begin{pmatrix} u_1, \dots, u_k \\ v_1, \dots, v_k \end{pmatrix} = \det(K(u_i, v_j))_{1 \leq i,j \leq k}.$$
Furthermore,
$$K\begin{pmatrix} u \\ v \end{pmatrix} = K\begin{pmatrix} u_1, \dots, u_k \\ v_1, \dots, v_k \end{pmatrix}.$$
Assume that $\psi$ is determinantal. For every $k \geq 1$ and every $z \in \mathbb{E}^k$ such that $K\binom{z}{z}$ is positive, let $\psi_z$ denote the determinantal process with kernel

$$K_z(u, v) = \frac{K\binom{u,z}{v,z}}{K\binom{z}{z}}, \qquad (u, v) \in \mathbb{E}^2. \tag{10}$$

With this notation, for every positive integer $k$ and $p$ and every $z \in \mathbb{E}^k$ and $v \in \mathbb{E}^p$,

$$K_z\begin{pmatrix} v \\ v \end{pmatrix} = \frac{K\binom{z,v}{z,v}}{K\binom{z}{z}}. \tag{11}$$

Note that, if $K\binom{z,v}{z,v}$ is positive, then $K\binom{z}{z}$ is positive and $\psi_{z,v} = (\psi_z)_v$.

A result of Shirai and Takahashi [27] (see also [19]) ensures that for $\mu_k$ almost every $z = (z_1, \dots, z_k) \in \mathbb{E}^k$,

$$\psi_z = (\psi \mid \{z_1, \dots, z_k\} \subset \psi) \setminus \{z_1, \dots, z_k\}. \tag{12}$$

Recall from [22], that if $\mathbb{E}$ is a vector space and $\psi$ is a process, stationary with respect to the translations of $\mathbb{E}$, then the associated Palm measure $\mathbb{Q}$ on $(\mathcal{M}_\sigma(\mathbb{E}), \mathcal{F})$ is defined by

$$\mathbb{Q}(M) = \frac{1}{\mu_1(A)} E \sum_{z \in \psi \cap A} 1_M \left( \sum_{z' \in \psi} \delta_{z'-z} \right), \qquad M \in \mathcal{F},$$

where $A \subset \mathbb{E}$ is an arbitrary Borel set such that $\mu_1(A)$ is positive and finite. It follows from (8) (see [27]) that the Palm measure of a determinantal stationary process $\psi$ with kernel $K$ is the law of the process $\psi_0 \cup \{0\}$ where $\psi_0 = (\psi \mid 0 \in \psi) \setminus \{0\}$ is determinantal with the kernel,

$$K_0(z_1, z_2) = \frac{K(z_1, z_2)K(0,0) - K(z_1, 0)K(0, z_2)}{K(0,0)}, \qquad (z_1, z_2) \in \mathbb{C}^2. \tag{13}$$

If $P\{\psi \neq \varnothing\}$ is positive, then $K(0,0)$ is positive. Applying this to the Ginibre kernel, one gets (2).

Note that if $\psi$ is a stationary Poisson process, that is, a point process with correlation functions $\rho_k \equiv 1/k!$ satisfying equality (9) for the degenerate, nonlocally trace-class kernel $K(z_1, z_2) = \delta_{z_1}(z_2)$, then applying formally the result by Shirai and Takahashi mentioned above, we obtain Slivnyak's formula [22], namely the fact that $\psi_{z_1,\ldots,z_k} \overset{\text{law}}{=} \psi$ for every positive $k$ and every distinct $z_j \in \mathbb{C}$.

For every Borel set $A$, let $N_A(\psi)$ denote the number of points of $\psi$ in $A$, that is,

$$N_A(\psi) = \sum_{z \in \psi} 1_A(z).$$

In the following, we assume that Conditions A and B below hold.

Condition A. The kernel $K(z_1, z_2)$ is a continuous function of $(z_1, z_2) \in \mathbb{E}^2$.

Condition B. For every bounded Borel set $A \subset \mathbb{E}$, the eigenvalues of the operator $K_A$ [acting on $L^2(A)$] are in the interval $[0, 1[$.

For every bounded Borel set $A \subset \mathbb{E}$, one sets

$$K_A^{(2)}(z_1, z_2) = \int_A K(z_1, v) K(v, z_2)\, d\lambda(v), \qquad (z_1, z_2) \in \mathbb{E}^2.$$

For every $n \geq 3$, $K_A^{(n)}$ denotes the iterated kernel of $K_A$, defined as

$$K_A^{(n)}(z_1, z_2) = \int_A K(z_1, v) K_A^{(n-1)}(v, z_2)\, d\lambda(v), \qquad (z_1, z_2) \in \mathbb{E}^2.$$

Conditions A and B above imply that the resolvent kernel,

$$R_A(z_1, z_2) = K(z_1, z_2) + \sum_{n \geq 2} K_A^{(n)}(z_1, z_2), \qquad (z_1, z_2) \in \mathbb{E}^2,$$

is a well-defined continuous function on $\mathbb{E}^2$.

REMARK 1. Note that the resolvent kernel is a continuous function of the domain in the sense that if $(A_n)_{n\geq 1}$ is a monotonous sequence of bounded Borel sets $A_n \subset E$ such that $A_n \uparrow A$ (and A is bounded) or $A_n \downarrow A$, then $R_{A_n}(z_1, z_2) \to R_A(z_1, z_2)$ $\lambda^2$ almost surely.

It is well known that, for every Borel set $A \subset E$, the probability of the event $\{N_A(\psi) = 0\}$ is a Fredholm determinant, namely,

$$P\{N_A(\psi) = 0\} = \det(I - K_A). \tag{14}$$

More generally Let $n \geq 1$ and $(A_i)_{1\leq i\leq n}$ denote $n$ disjoint, bounded Borel sets of positive measures $\lambda(A_i)$. Introduce

$$A = \bigcup_{i=1}^{n} A_i.$$

The Laplace transform of the joint law of random variables $N_{A_i}$, $i = 1, \ldots, n$, is given by the formula

$$E \exp\left(-\sum_{i=1}^{n} t_i N_{A_i}\right) = \det(I - K_{\bar{t},A}), \qquad t_i \in \mathbb{R}^+, i = 1, \ldots, n, \tag{15}$$

where $K_{\bar{t},A}$ designates the integral operator $K$ acting on the space $L^2(A, d\nu)$ with $d\nu(z) = \sum_{i=1}^{n}(1 - e^{-t_i})1_{A_i}(z)\, d\lambda(z)$. Now, (14) implies

$$P\{N_A(\psi) = 0\} = \exp\left\{-\int_A K(z,z)\, d\lambda(z) - \sum_{n\geq 2} \frac{1}{n} \int_A K_A^{(n)}(z,z)\, d\lambda(z)\right\} \tag{16}$$

and

$$P\{N_A(\psi) = 0\} = 1 + \sum_{n\geq 1} \frac{(-1)^n}{n!} \int_{A^n} K\begin{pmatrix} v \\ v \end{pmatrix} d\lambda^n(v). \tag{17}$$

The derivation of formulas (14)–(17) can be found in [27].

On the other hand, let us recall (see [25]) C. Platrier's classical formula from 1937 (established also by I. Fredholm for $k = 1$), that is, for every positive integer $k$ and every $z \in \mathbb{E}^k$, the relation,

$$\begin{aligned} K\begin{pmatrix} z \\ z \end{pmatrix} &+ \sum_{n\geq 1} \frac{(-1)^n}{n!} \int_{A^n} K\begin{pmatrix} z, v \\ z, v \end{pmatrix} d\lambda^n(v) \\ &= \left[1 + \sum_{n\geq 1} \frac{(-1)^n}{n!} \int_{A^n} K\begin{pmatrix} v \\ v \end{pmatrix} d\lambda^n(v)\right] R_A\begin{pmatrix} z \\ z \end{pmatrix}. \end{aligned} \tag{18}$$

From (11), (17) and (18), we deduce that, for every positive $k$ and every $z \in \mathbb{E}^k$ such that $K\binom{z}{z}$ is positive, for every bounded Borel set $A \subset \mathbb{E}$,

$$P\{N_A(\psi_z)=0\} = P\{N_A(\psi)=0\} \times \frac{R_A\binom{z}{z}}{K\binom{z}{z}}. \tag{19}$$

Remark 2. The kernel $R_A - K$ on $\mathbb{E}^2$ is obviously nonnegative. This fact, together with relation (19), implies that, for every positive integer $k$ and every $z \in \mathbb{E}^k$,

$$P\{N_A(\psi_z)=0\} \geq P\{N_A(\psi)=0\}. \tag{20}$$

We now establish some useful results. Let $n \geq 1$ and $(A_i)_{1\leq i\leq n}$ denote $n$ disjoint bounded Borel sets of positive measures $\lambda(A_i)$. Introduce

$$A = \bigcup_{i=1}^{n} A_i.$$

The following proposition gives the joint law of random variables $(N_{A_i})_{1\leq i\leq n}$.

Proposition 3. *Consider $n \geq 1$ nonnegative integers $k_i$ such that their sum $k = k_1 + \cdots + k_n$ is positive. Introduce*

$$B = \prod_{i=1}^{n} A_i^{k_i}, \qquad M = \{(N_{A_i}(\psi))_{1\leq i\leq n} = (k_i)_{1\leq i\leq n}\}.$$

*Then*

$$P\{M\} = \frac{P\{N_A(\psi)=0\}}{\prod_{i=1}^{n} k_i!} \int_B R_A\begin{pmatrix} z \\ z \end{pmatrix} d\lambda^k(z). \tag{21}$$

Proof. Observe that

$$P\{M\} = E \sum_{\tilde{z}\in\psi^{(k)}} f(\tilde{z},\psi), \tag{22}$$

where the function $f$ is defined as follows:

$$f(\tilde{z},\psi) = 1(N_A(\psi\setminus\tilde{z})=0) \sum_{(\tilde{z}_i)_i} \prod_{i=1}^{n} 1(\tilde{z}_i \subset A_i), \qquad |\tilde{z}| = k,$$

where the sums run above the following sets:

$$\bigcup_{i=1}^{n} \tilde{z}_i = \tilde{z} \qquad \forall 1 \leq i \leq n, |\tilde{z}_i| = k_i.$$

Now apply Campbell's formula (8). We obtain

$$
\begin{aligned}
P\{M\} &= E \sum_{\tilde{z}\in\psi^{(k)}} f(\tilde{z},\psi) \\
&= \int Ef(\{z_1,\dots,z_k\},(\psi \mid \{z_1,\dots,z_k\}\subset\psi))\,d\mu_k(z_1,\dots,z_k) \\
&= \int P\{N_A((\psi \mid \{z_1,\dots,z_k\}\subset\psi)\setminus\{z_1,\dots,z_k\})=0\} \\
&\qquad \times \sum_{(\tilde{z}_i)_i}\prod_{i=1}^{n} 1(\tilde{z}_i\subset A_i)\,d\mu_k(z_1,\dots,z_k).
\end{aligned}\tag{23}
$$

From property (12) we get

$$
\begin{aligned}
P\{M\} &= \int P\{N_A(\psi_z)=0\}\sum_{(\tilde{z}_i)_i}\prod_{i=1}^{n} 1(\tilde{z}_i\subset A_i)\,d\mu_k(z_1,\dots,z_k) \\
&= \frac{k!}{\prod_{i=1}^{n} k_i!}\int_B P\{N_A(\psi_z)=0\}\,d\mu_k(z_1,\dots,z_k),
\end{aligned}\tag{24}
$$

where $B=\prod_{i=1}^{n} A_i^{k_i}$.

The last equality above is obtained by counting partitions, noticing that $P\{N_A(\psi_z)=0\}$ depends on the set $\{z_1,\dots,z_k\}$ and that the measure $d\mu_k(z_1,\dots,z_k)$ is permutation invariant, that is, we have

$$
\int 1_D(z_1,\dots,z_k)\,d\mu_k(z_1,\dots,z_k)=\int 1_D(z_{\sigma(1)},\dots,z_{\sigma(k)})\,d\mu_k(z_1,\dots,z_k)
$$

for every Borel set $D\subset\mathbb{E}^k$ and for every permutation $\sigma\in\wp_k$.

Now, inserting formula (9) in (24) above we get

$$
P\{M\}=\frac{1}{\prod_{i=1}^{n} k_i!}\int_B P\{N_A(\psi_z)=0\}K\begin{pmatrix} z \\ z \end{pmatrix}\,d\lambda^n(z).
$$

It remains to apply formula (19) to obtain Proposition 3. □

Remark 4. Formula (23) works for any process $\psi$. In particular if we take for $\psi$ a Poisson process with intensity measure $\mu$ then, by Slivniak's formula, $\psi \overset{\text{law}}{=} (\psi \mid \{z_1,\dots,z_k\}\subset\psi)\setminus\{z_1,\dots,z_k\}$, the $k$th-order associated intensity measure is $d\mu_k(z_1,\dots,z_k)=(1/k!)\,d\mu(z_1)\cdots d\mu(z_k)$. Consequently, for a Poisson process $\psi$, formula (23) above gives the well-known expression

$$
\begin{aligned}
&P\{(N_{A_i}(\psi))_{1\le i\le n}=(k_i)_{1\le i\le n}\} \\
&\qquad = \frac{k!}{\prod k_i!}P\{N_A(\psi)=0\}\int_B d\mu_k(z_1,\dots,z_k)
\end{aligned}\tag{25}
$$

$$= \frac{\exp(-\mu(A))}{\prod k_i!} \prod_{i=1}^{n} \mu(A_i)^{k_i}.$$

Consider now $u \in \mathbb{E}$ such that $K(u,u)$ is positive. The process $\psi_u$ with kernel $K_u(z_1,z_2) = (1/K(u,u))[K(z_1,z_2)K(u,u) - K(z_1,u)K(u,z_2)]$, $(z_1,z_2) \in \mathbb{E}^2$, fulfills similar Conditions A and B. Indeed, if $K(z_1,z_2)$ is a continuous function of $(z_1,z_2) \in \mathbb{E}^2$, then $K_u(z_1,z_2)$ is a continuous function too. For every bounded Borel set $A \subset \mathbb{E}$, denote by $\alpha_{A,M}$ (resp. $\alpha_{u,A,M}$) the largest eigenvalue of the operator $K_A$ acting on $L^2(A)$ [resp. of the operator $K_{u,A}$ acting on $L^2(A)$]. Notice that the kernel

$$K(z_1,z_2) - K_u(z_1,z_2) = (1/K(u,u))K(z_1,u)K(u,z_2)$$

defines clearly a nonnegative operator, and thus $K \geq K_u$ in the Loewner order which implies inequality $\alpha_{u,A,M} \leq \alpha_{A,M}$. Consequently, if Condition B is satisfied by $K$ then it is satisfied by $K_u$ as well.

Denote by $R_{u,A}$ the associated resolvent kernel. Applying the relation (18) to the kernel $K_u$, one obtains that, for every $z \in \mathbb{E}^k$,

$$\begin{aligned} &K_u\begin{pmatrix} z \\ z \end{pmatrix} + \sum_{n\geq 1} \frac{(-1)^n}{n!} \int_{A^n} K_u\begin{pmatrix} z,v \\ z,v \end{pmatrix} d\lambda^n(v) \\ &\quad = \left[1 + \sum_{n\geq 1} \frac{(-1)^n}{n!} \int_{A^n} K_u\begin{pmatrix} v \\ v \end{pmatrix} d\lambda^n(v)\right] R_{u,A}\begin{pmatrix} z \\ z \end{pmatrix}. \end{aligned} \tag{26}$$

On the other hand,

$$\begin{aligned} &K\begin{pmatrix} u,z \\ u,z \end{pmatrix} + \sum_{n\geq 1} \frac{(-1)^n}{n!} \int_{A^n} K\begin{pmatrix} u,z,v \\ u,z,v \end{pmatrix} d\lambda^n(v) \\ &\quad = \left[1 + \sum_{n\geq 1} \frac{(-1)^n}{n!} \int_{A^n} K\begin{pmatrix} v \\ v \end{pmatrix} d\lambda^n(v)\right] R_A\begin{pmatrix} u,z \\ u,z \end{pmatrix}. \end{aligned} \tag{27}$$

From (10), (17), (26) and (27), we deduce that

$$R_A\begin{pmatrix} u,z \\ u,z \end{pmatrix} P\{N_A(\psi) = 0\} = K(u,u)R_{u,A}\begin{pmatrix} z \\ z \end{pmatrix} P\{N_A(\psi_u) = 0\}. \tag{28}$$

Applying formula (21) to the process $\psi_u$ and using (28) we obtain the proposition below.

PROPOSITION 5. *Consider $n$ nonnegative integers $k_i$ such that $k = k_1 + \cdots + k_n$ is positive. Introduce the set $B$ and the event $M_u$ defined as*

$$B = \prod_{i=1}^{n} A_i^{k_i}, \qquad M_u = \{(N_{A_i}(\psi_u))_{1\leq i\leq n} = (k_i)_{1\leq i\leq n}\}.$$

*Then,*

$$P\{M_u\} = \frac{P\{N_A(\psi)=0\}}{K(u,u)\prod_{i=1}^n k_i!}\int_B R_A\begin{pmatrix} u,z \\ u,z \end{pmatrix} d\lambda^k(z). \tag{29}$$

REMARK 6. Our notation for vectors of indices are such that equation (29) holds more generally, for every positive integer $p$ and every $u \in \mathbb{E}^p$ such that $K\binom{u}{u}$ is positive, if only one replaces the factor $K(u,u)$ in the denominator by $K\binom{u}{u}$.

We now state some simple consequences of Propositions 3 and 5. Denote, respectively, by $0 < \beta_n < 1$ and $h_n$, $n \geq 1$, the eigenvalues and the eigenfunctions of operator $K_A$. We recall that the eigenfunctions,

$$h_n(z) = \frac{1}{\beta_n}\int_A K(z,v)h_n(v)\,d\lambda(v), \qquad \|h_n\|_{L^2(A,d\lambda)} = 1, \tag{30}$$

are well defined and continuous on $\mathbb{E}$ and that

$$\det(I - K_A) = \prod_{n\geq 1}(1-\beta_n). \tag{31}$$

*In what follows, we shall always suppose that the eigenfunctions of operators are normalized* [*as in* (30)].

Assume now that $U \subset \mathbb{E}$ is a bounded open set and that $\lambda(V)$ is positive for every open subset $V \subset U$. Let $0 < \alpha_n < 1$ and $f_n$, $n \geq 1$, be the eigenvalues and the eigenfunctions of the operator $K_U$. A standard result (see, e.g., Theorem 2 of [30]) asserts the following.

LEMMA 7. *For every* $(z_1,z_2) \in \overline{U} \times \overline{U}$,

$$K(z_1,z_2) = \sum_{n\geq 1}\alpha_n f_n(z_1)\overline{f_n(z_2)}, \qquad R_U(z_1,z_2) = \sum_{n\geq 1}\frac{\alpha_n}{1-\alpha_n}f_n(z_1)\overline{f_n(z_2)},$$

*and the series are absolutely and uniformly convergent for* $z_1$ *and* $z_2$ *in every compact subset of* $U$.

REMARK 8. When the kernel $K$ has the form,

$$K(z_1,z_2) = \sum_{n=1}^M \beta_n h_n(z_1)\overline{h_n(z_2)}, \qquad (z_1,z_2) \in \mathbb{E}^2,$$

with functions $h_n$ that are continuous on $\mathbb{E}$ and orthonormal on a bounded Borel set $A$, then trivially,

$$R_A(z_1,z_2) = \sum_{n=1}^M \frac{\beta_n}{1-\beta_n}h_n(z_1)\overline{h_n(z_2)}, \qquad (z_1,z_2) \in \mathbb{E}^2.$$

Consider the case $n=1$ in Propositions 3 and 5. It was proved by Hough et al. ([12], Theorem 7) that the random variable $N_U(\psi)$ has the distribution of a sum of independent Bernoulli $(\alpha_i)$ random variables. Explicitly,

$$P\{N_U(\psi)=k\}=\sum_{(n_i)_i}\prod_{n\notin(n_i)_i}(1-\alpha_n)\prod_{i=1}^{k}\alpha_{n_i}, \tag{32}$$

where the sum runs over the indices $(n_i)_{1\le i\le k}$ such that $n_1<\cdots<n_k$. Now, by (14) and (31),

$$P\{N_U(\psi)=0\}=\det(I-K_U)=\prod_{n\ge 1}(1-\alpha_n)$$

and thus formula above can be written in the following form:

$$P\{N_U(\psi)=k\}=P\{N_U(\psi)=0\}\sum_{(n_i)_i}\prod_{1}^{k}\frac{\alpha_{n_i}}{1-\alpha_{n_i}}. \tag{33}$$

Assume that $u\in U$ and that $K(u,u)$ is positive. Using (29) we get

$$P\{N_U(\psi_u)=k\}=\frac{P\{N_U(\psi)=0\}}{K(u,u)}\Sigma_k, \tag{34}$$

where we introduce

$$\Sigma_k=\sum_{n\ge 1}|f_n(u)|^2\frac{\alpha_n}{1-\alpha_n}\sum_{(n_i)_i}\prod_{1}^{k}\frac{\alpha_{n_i}}{1-\alpha_{n_i}},$$

and where each last sum runs over the indices $(n_i)_{1\le i\le k}$ such that $n_1<\cdots<n_k$ and $n_i\ne n$ for every $1\le i\le k$.

Indeed, fix $M\ge 2$ and consider the kernels

$$K_{U,M}(z_1,z_2)=\sum_{n=1}^{M}\alpha_n f_n(z_1)\overline{f_n(z_2)}$$

and

$$R_{U,M}(z_1,z_2)=\sum_{n=1}^{M}\frac{\alpha_n}{1-\alpha_n}f_n(z_1)\overline{f_n(z_2)}, \qquad (z_1,z_2)\in U\times U,$$

where the functions $f_n$ are orthonormal on $U$. We have

$$\int_{U^k}R_{U,M}\begin{pmatrix}u,z\\u,z\end{pmatrix}d\lambda^k(z)$$

$$=\sum_{n=1}^{M}f_n(u)\frac{\alpha_n}{1-\alpha_n}\sum_{(n_i)_i}\prod_{1}^{k}\frac{\alpha_{n_i}}{1-\alpha_{n_i}} \tag{35}$$

$$\times \sum_{\sigma\in\wp_k} \int_{U^k} \prod_{j=1}^{k} f_{n_{\sigma(j)}}(z_j) \det \begin{pmatrix} \overline{f_n(u)} \cdots \overline{f_n(z_k)} \\ \overline{f_{n_{\sigma(1)}}(u)} \cdots \overline{f_{n_{\sigma(1)}}(z_k)} \\ \cdots \\ \overline{f_{n_{\sigma(k)}}(u)} \cdots \overline{f_{n_{\sigma(k)}}(z_k)} \end{pmatrix} d\lambda^k(z),$$

where sums run over the indices $(n_i)_{1\le i\le k}$ such that $1 \le n_1 < \cdots < n_k \le M$ and $n_i \ne n$ for every $1 \le i \le k$.

Observe that

$$\int_{U^k} \prod_{j=1}^{k} f_{n_{\sigma(j)}}(z_j) \det \begin{pmatrix} \overline{f_n(u)} \cdots \overline{f_n(z_k)} \\ \overline{f_{n_{\sigma(1)}}(u)} \cdots \overline{f_{n_{\sigma(1)}}(z_k)} \\ \cdots \\ \overline{f_{n_{\sigma(k)}}(u)} \cdots \overline{f_{n_{\sigma(k)}}(z_k)} \end{pmatrix} d\lambda^k(z) = \overline{f_n(u)}$$

due to the fact that the functions $f_n$ are orthonormal on $U$ (and $n_i \ne n$ for every $1 \le i \le k$).

Letting $M \to +\infty$ and applying Lemma 7, (29) and (35), we obtain formula (34).

Now, by elementary (but somewhat lengthy) computations, which we will not detail, we obtain the following proposition.

PROPOSITION 9. *With the assumptions above,*

$$P\{N_U(\psi_u) \le k\} - P\{N_U(\psi) \le k\} = \frac{P\{N_U(\psi) = 0\}}{K(u,u)} \tilde{\Sigma}_k,$$

*where*

$$\tilde{\Sigma}_k = \sum_{n\ge 1} |f_n(u)|^2 \frac{(\alpha_n)^2}{1-\alpha_n} \sum_{(n_i)_i} \prod_{1}^{k} \frac{\alpha_{n_i}}{1-\alpha_{n_i}},$$

*and where each last sum runs over the indices $(n_i)_{1\le i\le k}$ such that $n_1 < \cdots < n_k$ and $n_i \ne n$ for every $1 \le i \le k$.*

Proposition 9 implies the result below.

COROLLARY 10. *Let $U$ be a bounded open set and let $0 < \alpha_M < 1$ denote the largest eigenvalue of the operator $K_U$. Then for every nonnegative integer $k$, every positive integer $p$ and every $u \in U^p$ such that $K\binom{u}{u}$ is positive,*

$$P\{N_U(\psi_u) \le k\} \le (1-\alpha_M)^{-p} P\{N_U(\psi) \le k\}.$$

PROOF. By induction, if $u \in U$, then Proposition 9 and formula (34) imply

$$P\{N_U(\psi_u) \le k\} - P\{N_U(\psi) \le k\} \le \alpha_M P\{N_U(\psi_u) = k\}.$$

Therefore,

$$P\{N_U(\psi_u) \le k\} \le (1-\alpha_M)^{-1} P\{N_U(\psi) \le k\}.$$

Consider now $u = (v,w) \in U \times U^p$. Recall that $\psi_u = (\psi_w)_v$ and that $K \ge K_w$ in the Loewner order which implies inequality $\alpha_{w,M} \le \alpha_M$ where $\alpha_{w,M}$ denote the largest eigenvalue of the operator $K_{w,A}$. Thus

$$P\{N_U(\psi_u) \le k\} \le (1-\alpha_M)^{-1} P\{N_U(\psi_w) \le k\}$$

from which we obtain the announced result. □

REMARK 11. If $u \in \mathbb{E}$ and $z = (z_i)_{1\le i\le N} \in \mathbb{E}^N$, write $R_A^K\binom{u,z}{u,z}$ for the determinant $R_A\binom{u,z}{u,z}$ in which one replaces the terms $R_A(u,u)$ and $R_A(z_i,u)$ of the first column by $R_A(u,u) - K(u,u)$ and $R_A(z_i,u) - K(z_i,u)$, respectively. Then,

$$P\{N_A(\psi_u) \le k\} = P\{N_A(\psi) \le k\} + \frac{P\{N_A(\psi) = 0\}}{k!K(u,u)} \int_{A^k} R_A^K \begin{pmatrix} u, z \\ u, z \end{pmatrix} d\lambda^k(z).$$

One can prove this formula, from Propositions 3 and 5, by induction.

A further simple consequence of Proposition 3 and Remark 8 is the following. Consider another locally compact Hausdorff space $\mathbb{E}'$ with reference measure $\lambda'$, some bounded Borel nonintersecting sets $B_i \subset \mathbb{E}'$ of positive measures $\lambda'(B_i)$ and a point process $\psi' \subset \mathbb{E}'$ with kernel

$$L_B(z_1,z_2) = \sum_{n=1}^{M} \alpha_n g_n(z_1)\overline{g_n(z_2)}, \qquad (z_1,z_2) \in \mathbb{E}' \times \mathbb{E}',$$

where $0 < \alpha_n < 1$ and the functions $g_n$ for $1 \le n \le M$ are defined on $\mathbb{E}'$ and are orthonormal on

$$B = \bigcup_{i=1}^{N} B_i.$$

Now, let $\psi \subset \mathbb{E}$ be a point process with kernel

$$K_A(z_1,z_2) = \sum_{n=1}^{M} \alpha_n f_n(z_1)\overline{f_n(z_2)}, \qquad (z_1,z_2) \in \mathbb{E}^2,$$

where the functions $f_n$ for $1 \le n \le M$ are continuous on $\mathbb{E}$ and are orthonormal on

$$A = \bigcup_{i=1}^{N} A_i.$$

Assume that the following holds.

For every $1 \le i \le N$ and $1 \le n, m \le M$,

$$\int_{A_i} f_n(z)\overline{f_m(z)}\, d\lambda(z) = \int_{B_i} g_n(z)\overline{g_m(z)}\, d\lambda'(z). \tag{36}$$

Then the following proposition holds.

Proposition 12. *With the assumptions above, for every* $(k_i)_i$,

$$P\{(N_{A_i}(\psi))_{1\le i\le N} = (k_i)_{1\le i\le N}\} = P\{(N_{B_i}(\psi'))_{1\le i\le N} = (k_i)_{1\le i\le N}\}.$$

Proof. Let $k = k_1 + \cdots + k_N$. When $k = 0$, the result follows from formulas (14) and (31). Suppose now that $k$ is positive. By Remark 8 we have

$$R_A(z_1, z_2) = \sum_{n=1}^{M} \frac{\alpha_n}{1-\alpha_n} f_n(z_1)\overline{f_n(z_2)}, \qquad (z_1, z_2) \in A \times A,$$

and

$$R_B(z_1, z_2) = \sum_{n=1}^{M} \frac{\alpha_n}{1-\alpha_n} g_n(z_1)\overline{g_n(z_2)}, \qquad (z_1, z_2) \in B \times B.$$

Thus, for $\sigma \in \wp_k$ we obtain

$$\prod_{j=1}^{k} R_A(z_j, z_{\sigma(j)}) = \sum_{n_1,\dots,n_k=1}^{M} \prod_{j=1}^{k} \frac{\alpha_{n_j}}{1-\alpha_{n_j}} f_{n_j}(z_j)\overline{f_{n_{\sigma^{-1}(j)}}(z_j)}$$

and

$$\prod_{j=1}^{k} R_B(z_j, z_{\sigma(j)}) = \sum_{n_1,\dots,n_k=1}^{M} \prod_{j=1}^{k} \frac{\alpha_{n_j}}{1-\alpha_{n_j}} g_{n_j}(z_j)\overline{g_{n_{\sigma^{-1}(j)}}(z_j)}.$$

Denote $C = \prod_{i=1}^{n} A_i^{k_i}$ and $C' = \prod_{i=1}^{n} B_i^{k_i}$. Formula (36) implies that

$$\begin{aligned} &\int_C \prod_{j=1}^{k} f_{n_j}(z_j)\overline{f_{n_{\sigma^{-1}(j)}}(z_j)}\, d\lambda(z_1)\cdots d\lambda(z_k) \\ &\quad = \int_{C'} \prod_{j=1}^{k} g_{n_j}(z_j)\overline{g_{n_{\sigma^{-1}(j)}}(z_j)}\, d\lambda'(z_1)\cdots d\lambda'(z_k). \end{aligned} \tag{37}$$

Then, expanding the determinants appearing below and using the point (37) above, one gets the equality

$$\int_C R_A\binom{z}{z}\, d\lambda^k(z) = \int_{C'} R_B\binom{z}{z}\, d\lambda'^k(z). \tag{38}$$

This and (21) give the result. □

**3. Stochastic domination, proof of Theorem 3.** In this section, we assume that Condition I stated in the Introduction is satisfied.

Recall that a point process $\alpha \in \mathcal{M}_\sigma(\mathbb{E})$ stochastically dominates a point process $\beta \in \mathcal{M}_\sigma(\mathbb{E})$ if $Ef(\alpha) \geq Ef(\beta)$ for every bounded increasing measurable function $f$ defined on the space $(\mathcal{M}_\sigma(\mathbb{E}), \mathcal{F})$. It is is well known [17] that the point process $\alpha \in \mathcal{M}_\sigma(\mathbb{E})$ stochastically dominates the point process $\beta \in \mathcal{M}_\sigma(\mathbb{E})$ if and only if $P\{\alpha \in \mathcal{A}\} \leq P\{\beta \in \mathcal{A}\}$ for every decreasing event $\mathcal{A} \in \mathcal{F}$. Consider elementary decreasing events of the form $\{\forall 1 \leq i \leq M, N_{A_i} \leq k_i\} \in \mathcal{F}$ where $M$ is a positive integer, $k_i$, $1 \leq i \leq M$, are nonnegative integers and $A_i \subset \mathbb{E}$ are disjoint, bounded, Borel sets.

Denote by $\mathcal{F}_d \subset \mathcal{F}$ the collection of sets which are a finite union of such elementary decreasing events. The following lemma provides a useful tool in order to investigate stochastic domination properties of point processes.

Lemma 13. *The process $\beta$ is stochastically dominated by the process $\alpha$ if and only if, for every $\mathcal{A} \in \mathcal{F}_d$,*

$$P\{\alpha \in \mathcal{A}\} \leq P\{\beta \in \mathcal{A}\}.$$

Remark 14. The proof of Lemma 13 is standard. Similar characterizations are described, for example, in [17]; however, as pointed out by Yogeshwaran Dhandapani at ENS-DI-TREC (France), this result is not explicitly enunciaded in [17]. For completeness we sketch the proof of it in the Appendix.

We will now prove Theorem 3. Consider two kernels $K$ and $L$, satisfying Conditions A and B stated in the Introduction, such that $L \leq K$ in the Loewner order. Denote by $\varphi$ the process with kernels $L$ and by $\psi$ the process with kernel $K$. The idea of the proof is the following. By Lemma 13 we need to show that for every $\mathcal{A} \in \mathcal{F}_d$,

$$P\{\psi \in \mathcal{A}\} \leq P\{\varphi \in \mathcal{A}\}. \tag{39}$$

Fix the set $\mathcal{A} \in \mathcal{F}_d$. Applying the inclusion–exclusion principle it is easy to see that there exist nonintersecting, bounded Borel sets $B_i \subset \mathbb{E}$, $1 \leq i \leq N$, such that $P\{\psi \in \mathcal{A}\}$ (resp. $P\{\varphi \in \mathcal{A}\}$) can be expressed as a finite sum, up to the sign, of terms of the form,

$$P\{\forall i \in S, N_{B_i}(\psi) = k_i\} \qquad (\text{resp. } P\{\forall i \in S, N_{B_i}(\varphi) = k_i\}),$$

where $S \subset \{1, \ldots, N\}$. Let $U$ be an open bounded set such that

$$U \supset \bigcup_{i=1}^{N} B_i;$$

denote also, $B_0 = U \setminus \bigcup_{i=1}^N B_i$.

By Lemma 7 we have the spectral decomposition,

$$L_U(z_1, z_2) = \sum_{n \geq 1} \beta_n g_n(z_1) \overline{g_n(z_2)},$$

$$K_U(z_1, z_2) = \sum_{n \geq 1} \alpha_n f_n(z_1) \overline{f_n(z_2)}, \qquad (z_1, z_2) \in U.$$

The fact that $K \geq L$ in the Loewner order reads.

For all $f \in L^2(U)$,

$$\sum_{n \geq 1} \alpha_n \left| \int_U f_n(z) \overline{f(z)}\, d\lambda(z) \right|^2 \geq \sum_{n \geq 1} \beta_n \left| \int_U g_n(z) \overline{f(z)}\, d\lambda(z) \right|^2. \tag{40}$$

The inequality above implies that for every $n \geq 1$, the function $g_n$ is of the form $g_n = \sum_{k \geq 1} a_k^n f_k \in L^2(U)$. Denote $g_{n,M} = \sum_{k=1}^M a_k^n f_k$ and consider the nonnegative kernels,

$$K_{U,M} = \sum_{n=1}^M \alpha_n f_n(z_1) \overline{f_n(z_2)}$$

and

$$L_{U,M} = \sum_{n=1}^M \beta_n g_{n,M}(z_1) \overline{g_{n,M}(z_2)}, \qquad (z_1, z_2) \in U,$$

acting on $L^2(U)$.

Note that $\|K_{U,M}\| \leq \|K_U\| < 1$ and $\|L_{U,M}\| \leq \|L_U\| < 1$ where $\|\cdot\|$ denotes the supremum (operator) norm. Furthermore, if $\mathcal{V}^{(M)}$ is the subspace of $L^2(U)$ spanned by the functions $f_n$ with $1 \leq n \leq M$, then by (40), for each function $f \in \mathcal{V}^{(M)}$,

$$\begin{aligned} \sum_{n=1}^M \beta_n \left| \int_U g_{n,M}(z) \overline{f(z)}\, d\lambda(z) \right|^2 &= \sum_{n=1}^M \beta_n \left| \int_U g_n(z) \overline{f(z)}\, d\lambda(z) \right|^2 \\ &\leq \sum_{n=1}^M \alpha_n \left| \int_U f_n(z) \overline{f(z)}\, d\lambda(z) \right|^2. \end{aligned} \tag{41}$$

Denote by $\gamma_n$ and and $h_n$ the eigenvalues and the normalized eigenvectors of the operator $L_{U,M}$ [acting on $L^2(U)$]. The properties above imply that

$$0 \leq \gamma_n < 1 \quad \text{and} \quad h_n = \sum_{k=1}^M b_k^n f_k \in \mathcal{V}^{(M)}, \qquad 1 \leq n \leq M. \tag{42}$$

At last,

$$L_{U,M} = \sum_{n=1}^{M} \beta_n g_{n,M}(z_1)\overline{g_{n,M}(z_2)} = \sum_{n=1}^{M} \gamma_n h_n(z_1)\overline{h_n(z_2)}, \qquad (z_1, z_2) \in U.$$

Let $\varphi^{(M)} \subset U$ and $\psi^{(M)} \subset U$ be the processes associated, respectively, to the kernels $L_{U,M}$ and $K_{U,M}$.

LEMMA 15. *When $M \to \infty$,*

$$P\{\forall i \in S, N_{B_i}(\psi^{(M)}) = k_i\} \to P\{\forall i \in S, N_{B_i}(\psi) = k_i\},$$

$$P\{\forall i \in S, N_{B_i}(\varphi^{(M)}) = k_i\} \to P\{\forall i \in S, N_{B_i}(\varphi) = k_i\}.$$

PROOF. The straightforward consequence of (21), (31) and Lemma 7. □

It follows from Lemmas 13 and 15 that in order to prove that the process $\psi$ dominates the process $\varphi$ it suffices to show that, for every $M \geq 1$ inequality (39) is unchanged if we replace the terms of the form

$$P\{\forall i \in S, N_{B_i}(\psi) = k_i\} \qquad (\text{resp. } P\{\forall i \in S, N_{B_i}(\varphi) = k_i\}),$$

by the terms

$$P\{\forall i \in S, N_{B_i}(\psi^{(M)}) = k_i\} \qquad (\text{resp. } P\{\forall i \in S, N_{B_i}(\varphi^{(M)}) = k_i\}).$$

To obtain this result we use the fact that the stochastic domination occurs in the finite discrete determinantal process setting. See Theorem 6.2 and Paragraph 8 of [19], errata to [19] on Russel Lyon's website, and [2]. The link between our situation and a discrete determinantal process is given by the following lemma.

LEMMA 16. *Let $B_i$ denote nonintersecting Borel bounded subsets of $\mathbb{E}$, and let*

$$U = \bigcup_{i=0}^{N} B_i.$$

*Consider an orthonormal set of functions $\{l_n, n = 1, \ldots, M\} \subset L^2(U)$. Let $\mathcal{N}_i$ denote the dimension of the subspace $V_i \subset L^2(U)$ spanned by the functions $l_n 1_{B_i}$ with $1 \leq n \leq M$.*

*Then, there exists orthonormal vectors $z^n = (\mathrm{z}^n_{(0)}, \ldots, \mathrm{z}^n_{(N)})$, $z^n \in \prod_{i=0}^{N} \mathbb{C}^{\mathcal{N}_i}$, for $1 \leq n \leq M$, such that the following property holds.*

*For every $0 \leq i \leq N$ and every $1 \leq n, m \leq M$,*

$$\sum_{j=1}^{N_i} \mathrm{z}^n_{(i)}(j)\overline{\mathrm{z}^m_{(i)}(j)} = \int_{B_i} l_n(z)\overline{l_m(z)}\, dz. \tag{43}$$

PROOF. Since the sequence $(l_n)_n$ is orthonormal, property (43) implies that the sequence $(z^n)_n$ is orthonormal as well. Introduce an orthonormal basis $(e^i_j)_{1\le j\le \mathcal{N}_i}$ of the vector space $V_i \subset L^2(U)$. Then,

$$l_n 1_{B_i} = \sum_{j=1}^{\mathcal{N}_i} \lambda^n_{i,j} e^j_i.$$

The sequence defined by $\mathrm{z}^n_{(i)}(j) = \lambda^n_{i,j}$ fulfills property (43). $\square$

Consider now the vectors $z_n, n = 1, \ldots, N$, associated with Lemma 15, to the eigenvectors $f_n, n = 1, \ldots, M$, of the kernel $K_{U,M}$, and introduce the vectors,

$$v^n = (\mathrm{v}^n_{(0)}, \ldots, \mathrm{v}^n_{(N)}) = \sum_{k=1}^{M} b^n_k z^k \in \prod_{i=0}^{N} \mathbb{C}^{\mathcal{N}_i}, \qquad n = 1, \ldots, M,$$

related to functions $h_n, n = 1, \ldots, M$, of (42). Notice that for every $0 \le i \le N$ and every $1 \le n, m \le M$,

$$\sum_{j=1}^{\mathcal{N}_i} \mathrm{v}^n_{(i)}(j)\overline{\mathrm{v}^m_{(i)}(j)} = \int_{B_i} h_n(z)\overline{h_m(z)}\, d\lambda(z). \tag{44}$$

Consequently, $(v^n)_n$ is a set of orthonormal vectors.

Consider now the sets

$$E_j = \{(i,j); 1 \le i \le \mathcal{N}_j\}, \qquad 0 \le j \le N, E = \bigcup_{j=0}^{N} E_j,$$

and the discrete kernels defined on $E$ by,

$$K_M((i_1, j_1), (i_2, j_2)) = \sum_{n=1}^{M} \alpha_n \mathrm{z}^n_{(i_1)}(j_1)\overline{\mathrm{z}^n_{(i_2)}(j_2)}$$

and

$$L_M((i_1, j_1), (i_2, j_2)) = \sum_{n=1}^{M} \gamma_n \mathrm{v}^n_{(i_1)}(j_1)\overline{\mathrm{v}^n_{(i_2)}(j_2)}.$$

Inequality (41) implies that $L_M \le K_M$ in the Loewner order. Introduce the determinantal process $\chi \subset E$ with kernel $L_M$, and the process $\zeta \subset E$ with kernel $K_M$. Proposition 12 and formulas (43), (44) imply that

$$P\{\forall i \in S, N_{B_i}(\psi^{(M)}) = k_i\} = P\{\forall i \in S, N_{E_i}(\zeta) = k_i\}$$

and

$$P\{\forall i \in S, N_{B_i}(\varphi^{(M)}) = k_i\} = P\{\forall i \in S, N_{E_i}(\chi) = k_i\}.$$

Consequently if we replace in formula (39), the terms of the form

$$P\{\forall i \in S, N_{B_i}(\varphi) = k_i\} \qquad (\text{resp. } P\{\forall i \in S, N_{B_i}(\psi) = k_i\})$$

by the terms

$$P\{\forall i \in S, N_{B_i}(\varphi^{(M)}) = k_i\} \qquad (\text{resp. } P\{\forall i \in S, N_{B_i}(\psi^{(M)}) = k_i\}),$$

we obtain inequality

$$P\{\zeta \in \mathcal{A}'\} \le P\{\chi \in \mathcal{A}'\}$$

for a suitable decreasing event $\mathcal{A}' \in \mathcal{F}_d(E)$. The above mentioned result of [2] and [19] asserts that this inequality is indeed true and thus the proof of Theorem 3 is finished.

We are now in position to apply the celebrated Strassen's theorem. This follows from the fact that the space $\mathcal{M}_\sigma(\mathbb{E})$ of counting measures endowed with the vague topology is a Polish space and its associated Borel $\sigma$-algebra coincides precisely with the $\sigma$-algebra $\mathcal{F}$ (see [17] and [18]).

THEOREM 6. *With the hypothesis of Theorem 3, there exists a point process $\eta$ such that*

$$\psi \overset{\text{law}}{=} \varphi \cup \eta, \qquad \varphi \cap \eta = \varnothing.$$

Theorem 3 implies Theorem 2. More generally we have:

THEOREM 7. *Let $\psi$ be a point process satisfying Conditions A and B of the Introduction. For all points $u$ such that $K(u,u)$ is positive, the process $\psi$ dominates stochastically the process $\psi_u$.*

PROOF. It is obvious that $K \ge K_u$ in the Loewner order. □

PROBLEM 1. Prove Theorem 7 directly from (21) and (29).

**4. Palm measure of the Ginibre process, proof of Theorem 1.** Recall that the Ginibre process $\phi \subset \mathbb{R}^2 = \mathbb{C}$ is a stationary, isotropic point process satisfying Conditions A and B of the Introduction. The reference measure $\lambda$ is the area measure of $\mathbb{R}^2$, and Condition I is trivially satisfied. Moreover, for every integer $k \ge 1$ and every set of distinct points $\{z_1, \dots, z_k\}$ included in $\mathbb{C}$, respectively, in $\mathbb{C} \setminus \{0\}$, $K\binom{z_1,\dots,z_k}{z_1,\dots,z_k}$ is positive, respectively, $K_0\binom{z_1,\dots,z_k}{z_1,\dots,z_k}$ is positive.

From formula (13) it follows that the process $\phi_0 = (\phi \mid 0 \in \phi) \setminus \{0\}$ is determinantal with the kernel $K_0$ such that

$$K_0(z_1, z_2) = (1/\pi)(e^{z_1 \overline{z}_2} - 1)\exp(-(1/2)(|z_1|^2 + |z_2|^2)), \qquad (z_1, z_2) \in \mathbb{C}^2.$$

The intensity measure $\mu_{0,1}$ is absolutely continuous with respect to the measure $\lambda$ and has for density the correlation function

$$K_0(z,z) = (1/\pi)(1 - e^{-|z|^2}). \tag{45}$$

In particular, the process $\phi_0$ is not stationary.

Remark 17. The stationarity of the Ginibre process $\phi$ is expressed by the fact that for each fixed $a \in \mathbb{C}$ the determinantal point process with kernel $\widehat{K}$ such that $\widehat{K}(z_1, z_2) = K(z_1 - a, z_2 - a)$, that is,

$$\widehat{K}(z_1, z_2) = (1/\pi) e^{(z_1-a)(\overline{z_2-a}) - (1/2)(|z_1-a|^2+|z_2-a|^2)}, \qquad (z_1, z_2) \in \mathbb{C}^2,$$

coincides (in law) with $\phi$. Note that $K \neq \widehat{K}$.

Consider now the radial processes $|\phi|$ and $|\phi_0|$. The result below is well known [12, 15].

Theorem 8 (Kostlan). *Let $X_{n,m}$ with $n \geq 1$ and $m \geq 1$ denote i.i.d. random variables with exponential distribution $e^{-x}\,dx$ on $x \geq 0$. For every $n \geq 1$, let*

$$R_n = \sqrt{\sum_{m=1}^{n} X_{n,m}}.$$

*Then, the collection of moduli of the points of $\phi$ has the same distribution as the collection of random variables $\{R_n, n \geq 1\}$.*

$$|\phi| \stackrel{\text{law}}{=} \{R_n, n \geq 1\}. \tag{46}$$

Remark 18. Note that Theorem 8 implies that, almost surely, there exists no $(z_1, z_2) \in \phi \times \phi$ such that $z_1 \neq z_2$ and $|z_1| = |z_2|$.

We will shown that result (46) can be deduced from formula (15). Indeed, let us fix $0 < r_1 < \cdots < r_n = r$ and consider the sets $A_1 = B(r_1)$, $A_i = \{z \in \mathbb{C}; r_{i-1} < |z| \leq r_i\}$ for $i = 2, \ldots, n$, $B(r) = \bigcup_{i=1}^{n} A_i$. Also, let us fix $t_i > 0$, $i = 1, \ldots, n$.

Observe that the functions,

$$f_n(z) = (1/\sqrt{\pi n!}) e^{-(1/2)|z|^2} z^n, \qquad z \in B(r), n \geq 1,$$

are orthogonal on $B(r)$ with respect to the measure,

$$d\nu(z) = \sum_{i=1}^{n} (1 - e^{-t_i}) 1_{A_i}(z)\, d\lambda(z).$$

Denote $\alpha_n = \int_{B(r)} |f_n|^2 \, d\nu(z)$ then normalizing, we obtain

$$K(z_1, z_2) = (1/\pi) e^{z_1 \overline{z}_2 - (1/2)(|z_1|^2 + |z_2|^2)} = \sum_{n \geq 1} \alpha_n \widehat{f_n}(z_1) \overline{\widehat{f_n}(z_2)}$$

with $\widehat{f_n}(z) = (1/\sqrt{\alpha_n}) f_n$. Consider now the radial process $|\phi|$ and the intervals $I_1 = [0, r_1], \ldots, I_n = ]r_{n-1}, r_n]$. Formulas (15) and (31) imply that

$$\begin{aligned} E \exp\left(-\sum_{i=1}^{n} t_i N_{I_i}(|\phi|)\right) &= E \exp\left(-\sum_{i=1}^{n} t_i N_{A_i}(\phi)\right) \\ &= \det(I - K_{\bar{t}, A}) = \prod_{n \geq 1} (1 - \alpha_n). \end{aligned} \tag{47}$$

Computing (what is an elementary exercise) the Laplace transform $E \exp(-\sum_{i=1}^{n} t_i N_{I_i}(\mathcal{R}))$ for the point process $\mathcal{R} = \{R_n, n \geq 1\}$ gives exactly the same value. Thus $|\phi| \stackrel{\text{law}}{=} \{R_n, n \geq 1\}$. More generally, if $\psi(F)$, $F \subset \mathbb{N}$, $F \neq \varnothing$, is the point process related to the kernel

$$K(F)(z_1, z_2) = \sum_{n \in F} \alpha_n \widehat{f_n}(z_1) \overline{\widehat{f_n}(z_2)},$$

then

$$|\phi(F)| \stackrel{\text{law}}{=} \{R_n, n \in F\}. \tag{48}$$

In particular,

$$|\phi_0| \stackrel{\text{law}}{=} \{R_n, n \geq 2\}$$

and

$$\{R_n, n \geq 1\} = \{R_n, n \geq 2\} \cup \{R_1\}$$

provide a disjoint coupling of $|\phi_0|$ and $\{R_1\}$ with union marginal $|\phi|$.

Consider now, for $M \geq 1$, the kernels

1. $K_M(z_1, z_2) = (1/\pi) \sum_{n=0}^{M} ((z_1 \overline{z}_2)^n / n!) \exp(-(1/2)(|z_1|^2 + |z_2|^2))$,
2. $K_{0,M}(z_1, z_2) = (1/\pi) \sum_{n=1}^{M} ((z_1 \overline{z}_2)^n / n!) \exp(-(1/2)(|z_1|^2 + |z_2|^2))$,

and denote by $\phi^{(M)}$ the point process associated with the kernel $K_M$ and by $\phi_0^{(M)}$ the point process associated with the kernel $K_{0,M}$.

Observe that on one hand we have

$$E \operatorname{card}\{\phi^{(M)}\} = \int_{\mathbb{C}} K_M(z, z) \, dz = M,$$

$$E \operatorname{card}\{\phi_0^{(M)}\} = \int_{\mathbb{C}} K_{0,M}(z, z) \, dz = M - 1,$$

and on the other hand, by (9), the correlation function of order $M+1$ (resp. $M$) for the process $\phi^{(M)}$ (resp. $\phi_0^{(M)}$) is equal to zero which implies that $\mathrm{card}\{\phi^{(M)}\} \le M$ and $\mathrm{card}\{\phi_0^{(M)}\} \le M-1$, almost surely. Therefore, $\mathrm{card}\{\phi^{(M)}\} = M$ and $\mathrm{card}\{\phi_0^{(M)}\} = M-1$, almost surely. Moreover, we have $K_{0,M} \le K_M$ in the Loewner order. Formula (48) implies also

$$(49)\quad |\phi^{(M)}| \overset{\text{law}}{=} \{R_n, 1 \le n \le M+1\}, \qquad |\phi_0^M| \overset{\text{law}}{=} \{R_n, 2 \le n \le M+1\}.$$

It follows from the properties above and from Theorem 6 that there exists a disjoint coupling,

$$(50)\qquad \phi^{(M)} = \phi_0^{(M)} \cup \eta^{(M)}, \qquad \phi_0^{(M)} \cap \eta^{(M)} = \varnothing,$$

such that the point process $\eta^{(M)}$ is a single random variable, $\eta^{(M)} = \{Z_M\}$.

By equation (49) and the fact that Remark 18 also applies to the process $\phi^{(M)}$, we deduce that $|Z_M| \overset{\text{law}}{=} R_1$.

Lemma 19. *The random variable $Z_M$ is centered Gaussian, and*

$$E|Z_M|^2 = 1.$$

Proof. The random variable $|Z_M|^2$ has exponential distribution $e^{-x}\,dx$, $x \ge 0$, thus it suffices to show that the law of $Z_M$ is invariant by rotations $O$ with its center at the origin, that is, $P\{Z_M \in A\} = P\{Z_M \in O(A)\}$, for every such $O$. Simple computation gives

$$(51)\quad P\{Z_M \in A\} = \sum_{0 \le k \le M+1} [P\{N_A(\phi_0^{(M)}) \le k\} - P\{N_A(\phi^{(M)}) \le k\}].$$

The processes $\phi^{(M)}$ and $\phi_0^{(M)}$ are isotropic, hence formula (51) implies the result.

Consider now the laws $P^{(M)}$, $M \ge 1$, of random elements $(Z_M, \phi_0^{(M)})$ with values in the product space $\mathbb{C} \times \mathcal{M}_\sigma(\mathbb{C})$ [the space $\mathcal{M}_\sigma(\mathbb{C})$ being endowed with the vague topology].

Denote, respectively, by $Q$, $Q^{(M)}$, $Q_0$ and $Q_0^{(M)}$, the laws of the processes $\phi$, $\phi^{(M)}$, $\phi_0$ and $\phi_0^{(M)}$. Finally, let $I : \mathbb{C} \times \mathcal{M}_\sigma(\mathbb{C}) \longrightarrow \mathcal{M}_\sigma(\mathbb{C})$ be the continuous application defined by $I(x,\zeta) = \{x\} \cup \zeta$. □

Lemma 20. *The following properties hold:*

1. *The sequences $(Q^{(M)})_M$ and $(Q_0^{(M)})_M$ are tight.*
2. $Q^{(M)} \xrightarrow[M\to+\infty]{D} Q$ *and* $Q_0^{(M)} \xrightarrow[M\to+\infty]{D} Q_0$.
3. *The sequence $(P^{(M)})_M$ is tight.*

4. *Consider the probability* $P^{(M)}$ *on* $(\mathbb{C}\times\mathcal{M}_\sigma(\mathbb{C}),\mathcal{B}(\mathbb{C})\otimes\mathcal{F})$, *then* $I\stackrel{D}{=}Q^{(M)}$.

PROOF. Property 1 is obvious from the characterization of tightness for random measures (see [14], page 33). Property 2 follows from property 1 and Proposition 3. Property 3 is a consequence of property 1 and the fact that, by Lemma 19, the standard normal law coincides with the marginal law on $\mathbb{C}$ of the probability $P^{(M)}$. Finally, property 4 is nothing but the equality in law (50). □

It is well known that a suitable subsequence of $(P^{(M)})_M$ converges in distribution to a probability $P^*$ on $\mathbb{C}\times\mathcal{M}_\sigma(\mathbb{C})$. Lemma 20 implies that $P^*$ has, for marginal laws, the standard normal law and $Q_0$ and that with $P^*$ on $\mathbb{C}\times\mathcal{M}_\sigma(\mathbb{C})$, we obtain $I\stackrel{D}{=}Q$. Consequently, a random element with distribution $P^*$ provides a disjoint coupling $(Z,\phi_0)$ of $\phi$. The proof of Theorem 1 is then finished.

One can notice also that we have

$$P\{Z\in A\}=\sum_{k\geq 0}[P\{N_A(\phi_0)\leq k\}-P\{N_A(\phi)\leq k\}];\tag{52}$$

thus, if $U$ is an open set containing the origin, then by inserting the formula of Proposition 9 in (52), we obtain

$$P\{Z\in U\}=\frac{P\{N_U(\phi)=0\}}{K(0,0)}\sum_{n\geq 1}|f_n(0)|^2\frac{(\alpha_n)^2}{1-\alpha_n}\left[1+\sum_{k\geq 1}\sum_{(n_i)_i}\prod_1^k\frac{\alpha_{n_i}}{1-\alpha_{n_i}}\right],$$

where the last sum is over the integers $(n_i)_{1\leq i\leq k}$ such that $n_i<n_{i+1}$ and $n_i\neq n$ for every $i$. Hence,

$$P\{Z\in U\}=\frac{P\{N_U(\phi)=0\}}{K(0,0)}\times\sum_{n\geq 1}(\alpha_n)^2|f_n(0)|^2\times\prod_{i\geq 1}\frac{1}{1-\alpha_i},$$

and, finally,

$$P\{Z\in U\}=\frac{K_U^{(2)}(0,0)}{K(0,0)}=\frac{1}{\pi}\int_U e^{-|z|^2}\,dz.$$

Thus we find again that the law of $Z$ is Gaussian. Notice also the formula,

$$P\{Z\in A\mid N_A(\phi_0)=0\}=1-\frac{K(0,0)}{R_A(0,0)},\tag{53}$$

which follows from (19) via the identities,

$$P\{Z\in A,N_A(\phi_0)=0\}=P\{N_A(\phi_0)=0\}-P\{N_A(\phi)=0\}$$

and

$$P\{N_A(\phi)=0\} = \frac{K(0,0)}{R_A(0,0)} P\{N_A(\phi_0)=0\}.$$

PROBLEM 2. This is an open problem, that is, to know how the random variable $Z$ is correlated with the point process $\phi_0$. A similar unsolved problem arises in the framework of finite discrete determinantal processes (see question (10.1) in [19]).

REMARK 21. The method we used to prove Theorem 1, that is, a coupling result [formula (50)] in a finite-dimensional case associated with a "tightness argument" (Lemma 20) is very similar to that used by R. Lyons, in the discrete determinantal process setting, to prove Proposition 10.3 in [19].

REMARK 22. Theorem 6 can be applied to the processes $\phi$ and $\phi_0$ and thus provides a disjoint coupling $\phi = \psi_0 \cup \{\eta\}$. However, there is a difficulty to deduce Theorem 1 directly from this (due to the fact that it is unclear that the process $\eta$ could be taken as being a single random variable).

REMARK 23. The random variables $R_n^2$, $n \geq 1$, are Gamma$(n,1)$ distributed and independent. They are stochastically increasing but not almost surely increasing. It is interesting to note that if $\widetilde{R}_n = \sqrt{\sum_{m=1}^{n} X_{1,m}}$, $n \geq 1$, is the radial process of a Poisson stationary process which has the same intensity $(1/\pi)\,dz$ as the process $\phi$, then the random variables $\widetilde{R}_n^2$, $n \geq 1$, are Gamma$(n,1)$ distributed as well; they are almost surely increasing and (of course) is not independent.

PROBLEM 3. Construct explicitly random variables $Z_n$, $n \geq 1$, such that:

1. $\phi \overset{\text{law}}{=} \{Z_n, n \geq 1\}$;
2. $\forall n \geq 1, |Z_n| \overset{\text{law}}{=} R_n$;
3. $\phi_0 \overset{\text{law}}{=} \{Z_n, n \geq 2\}$.

REMARK 24. The Palm measure of $\phi^{\star\alpha}$ is obtained by adding the origin and deleting the point $\sqrt{\alpha}Z$ if the latest belongs (which occurs with probability $\alpha$) to the process $\phi^{\star\alpha}$.

REMARK 25. Similar results could be proved for the point process in the unit disk of $\mathbb{C}$ related to the Bergman kernel and studied in [24].

**5. Ginibre–Voronoi tessellation, proof of Theorems 3 and 4.** Consider now the space $\mathcal{K}$ of compact convex sets of $\mathbb{R}^2 = \mathbb{C}$ endowed with the usual Hausdorff metric. For every point process $\psi$, let

$$C(u,\psi) = \{z \in \mathbb{C}; \forall v \in \psi, |z-u| \leq |z-v|\},$$

and let $\{C(u,\phi); u \in \phi\}$ denote the Voronoi tessellation generated by the Ginibre process $\phi$. Recall that its statistical properties, namely its empirical distributions (the process being ergodic), are described [3, 5, 6] by the typical cell $\mathcal{C}$ defined by means of the identity,

$$Eh(\mathcal{C}) = \frac{\pi}{\lambda(B)} E \sum_{z \in B \cap \psi} h(C(z) - z),$$

where $h$ runs through the space of positive measurable functions on $\mathcal{K}$, and $B \subset \mathbb{C}$ is an arbitrary Borel set with the finite positive area $\lambda(B)$. Consider now the cell

$$C(0,\phi_0) = \{z \in \mathbb{C}; \forall u \in \phi_0, |z| \leq |z-u|\}.$$

Campbell's formula (8) gives the identity,

$$Eh(\mathcal{C}) = \frac{\pi}{\lambda(B)} E \sum_{z \in B \cap \psi} h(C(0, \phi - z)) = Eh(C(0,\phi_0)).$$

Hence,

$$\mathcal{C} \overset{\text{law}}{=} C(0,\phi_0). \tag{54}$$

In what follows, we shall use the notation $C(0,\phi_0) = C(0)$. The law of the random set $C(0)$ can be obtained by means of the method described in [4]. Let us introduce some notation. Fix $k \geq 1$.

- For every $u \in \mathbb{C}$, let $H(u) = \{z \in \mathbb{C}; \langle z-u, u\rangle \leq 0\}$.
- For every $z \in \mathbb{C}^k$ with $z = (z_1, \dots, z_k)$, let $\mathfrak{H}(z)$ denote the intersection of half-spaces,

$$\mathfrak{H}(z) = \bigcap_{i=1}^{k} H(z_i/2).$$

- For every $z \in \mathbb{C}^k$, let $\mathcal{F}(z) = \bigcup_{u \in \mathfrak{H}(z)} B(u, |u|)$ where $B(u,r)$ denotes the disk centered at $u$ and of radius $r \geq 0$.
- Let $A \subset \mathbb{C}^k$ denote the set of $z \in \mathbb{C}^k$ such that $\mathfrak{H}(z)$ is a bounded polygon with $k$ sides.

THEOREM 9. *For every $k \geq 3$,*

$$P\{\mathcal{C} \text{ has } k \text{ sides}\} = \frac{1}{k!} \int_A P\{N_{\mathcal{F}(z)}(\phi_0) = 0\} R_{0,\mathcal{F}(z)} \binom{z}{z} \, dz.$$

PROOF. Observe that

$$P\{C(0) \text{ has } k \text{ sides}\} = E \sum_{\tilde{z} \in \phi_0^{(k)}} 1_A(\tilde{z}) \times 1_{\{0\}}(N_{\mathcal{F}(z_1,\ldots,z_k)}(\phi_0 \setminus \tilde{z})).$$

With Campbell's formula (8) applied to the process $\phi_0$, one gets

$$P\{\mathcal{C} \text{ has } k \text{ sides}\} = P\{C(0) \text{ has } k \text{ sides}\} = \int_A P\{N_{\mathcal{F}(z)}(\phi_{0,z}) = 0\}\, d\mu_{0,k}(z)$$

and hence by formulas, (9) and (19), we obtain

$$\begin{aligned} P\{C(0) \text{ has } k \text{ sides}\} &= \frac{1}{k!} \int_A P\{N_{\mathcal{F}(z)}(\phi_{0,z}) = 0\} K_0 \begin{pmatrix} z \\ z \end{pmatrix} dz \\ &= \frac{1}{k!} \int_A P\{N_{\mathcal{F}(z)}(\phi_0) = 0\} R_{0,\mathcal{F}(z)} \begin{pmatrix} z \\ z \end{pmatrix} dz. \end{aligned}$$

□

In the same way, one can compute, conditionally on the fact that the cell has $k$ sides, the expectation of an arbitrary, measurable, positive functional of $\mathcal{C}$ which is expressed through a function $f$ acting on points $\{z_1, \ldots, z_k\} = \mathcal{N}(\phi_0) \subset \phi_0$ for which the bisecting line of the interval $[0, z_i]$ intersects the cell $C(0)$. The resulting integral will have the form

$$\frac{1}{k!} \int_A f(z) P\{N_{\mathcal{F}(z)}(\phi_0) = 0\} R_{0,\mathcal{F}(z)} \begin{pmatrix} z \\ z \end{pmatrix} dz. \tag{55}$$

Deconditioning, one can obtain analytical formulas of the laws of the geometric characteristics of the typical cell $\mathcal{C}$. Note that formula (17) gives an analytical expression of the probability $P\{N_{\mathcal{F}(z)}(\phi_0) = 0\}$ which appears in (55). Unfortunately, these integrals are complicated and numerical computations are difficult. This drawback appears already in [4] for the typical cell of the Poisson–Voronoi tessellation.

A general result asserts [22] that the first-order moment of the area $V(\mathcal{C})$ of the cell $\mathcal{C}$ is equal to $EV(\mathcal{C}) = \pi$. The moments $EV^k(\mathcal{C})$ of higher orders can be expressed in terms of integrals more tractable than (55). Recall our notation

$$D(z_1, \ldots, z_k) = \bigcup_{i=1}^{k} B(z_i, |z_i|) \subset \mathbb{C}.$$

We use the fact that

$$z \in C(0)^k \quad \Longleftrightarrow \quad N_{D(z)}(\phi_0) = 0. \tag{56}$$

Let $A \subset \mathbb{C}$ be a Borel set. From (54), (56) and (19),

$$\begin{aligned} E[V^k(\mathcal{C} \cap A)] &= \int_{A^k} P\{N_{D(z)}(\phi_0) = 0\}\, dz \\ &= \int_{A^k} \frac{R_{D(z)}(0,0)}{K(0,0)} P\{N_{D(z)}(\phi) = 0\}\, dz \end{aligned} \tag{57}$$

and hence by (16), we obtain

$$E[V^k(\mathcal{C}\cap A)] = \int_{A^k} \exp\biggl\{-\int_{D(z)} K_0(u,u)\,du \qquad (58)$$

$$- \sum_{n\geq 2}\frac{1}{n}\int_{D(z)} K_{0,D(z)}^{(n)}(u,u)\,du\biggr\}\,dz$$

and

$$E[V^k(\mathcal{C}\cap A)] = \int_{A^k} \frac{R_{D(z)}(0,0)}{K(0,0)} \exp\biggl\{-\frac{1}{\pi}V(z) \qquad (59)$$

$$- \sum_{n\geq 2}\frac{1}{n}\int_{D(z)} K_{D(z)}^{(n)}(u,u)\,du\biggr\}\,dz,$$

where $V(z)$ denote the area of the set $D(z) = D(z_1,\ldots,z_k)$.

We will now use formulas (58) and (59) in order to compare the area of $\mathcal{C}$ with the area of the typical cell $\mathcal{C}_p$ of the Voronoi tessellation associated with a stationary Poisson process which has the same intensity measure $(1/\pi)\,dz$ as the process $\phi$. For every Borel set A,

$$E[V^k(\mathcal{C}_p\cap A)] = \int_{A^k} e^{-V(z)/\pi}\,dz. \qquad (60)$$

Hence by (58), (45) and the fact that for every $z = (z_1,\ldots,z_k) \in B(r)^k$, $D(z) = D(z_1,\ldots,z_k) \subset B(2r)$, one has

$$E[V^k(\mathcal{C}\cap B(r))] \leq E[V^k(\mathcal{C}_p\cap B(r))]\exp\biggl(\frac{1}{\pi}\int_{B(2r)} e^{-|z|^2}\,dz\biggr). \qquad (61)$$

5.1. *Proof of Theorem 4, part* (a). Formula (60) implies that

$$E[V^k(\mathcal{C}_p\cap B(r))] = \int_{B(r)^k} e^{-V(z)/\pi}\,dz = r^{2k}\int_{B(1)^k} e^{-r^2V(z)/\pi}\,dz. \qquad (62)$$

For $z\in\mathbb{C}^k$, let $\alpha_{0,n}$, $n\geq 1$, denote the eigenvalues of $K_0$ acting on $D(z)$ where $\alpha_{0,1}$ is the largest eigenvalue, then

$$\sum_{m\geq 2}\frac{1}{m}\int_{D(z)} K_{0,D(z)}^{(m)}(u,u)\,du = \sum_{m\geq 2}\frac{1}{m}\sum_{n\geq 1}(\alpha_{0,n})^m \qquad (63)$$

$$\leq \frac{1}{2(1-\alpha_{0,1})}\sum_{n\geq 1}(\alpha_{0,n})^2.$$

Furthermore,

$$\sum_{n\geq 1}(\alpha_{0,n})^2 = \int_{D(z)} K_{0,D(z)}^{(2)}(u,u)\,du = (1/\pi^2)\int_{D(z)^2} |1-e^{u\overline{v}}|^2 e^{-|u|^2-|v|^2}\,du\,dv.$$

If $z \in B(r)^k$, then $D(z) \subset B(2r)$, and hence

$$\sum_{n\geq 1}(\alpha_{0,n})^2 \leq (1/\pi^2)\int_{B(2r)^2} |1-e^{u\overline{v}}|^2 e^{-|u|^2-|v|^2}\,du\,dv = O(r^8).$$

This and (63) imply (since $\alpha_{0,1}$ decrease when domain decrease)

$$\sum_{m\geq 2}\frac{1}{m}\int_{D(z)} K^{(m)}_{0,D(z)}(u,u)\,du = O(r^8). \tag{64}$$

Therefore, by (58) we obtain

$$\begin{aligned} E[V^k(\mathcal{C}\cap B(r))] &= \int_{B(r)^k}\exp\Big\{-\int_{D(z)} K_0(u,u)\,du \\ &\qquad - \sum_{m\geq 2}\frac{1}{m}\int_{D(z)} K^{(m)}_{0,D(z)}(u,u)\,du\Big\}\,dz \\ &= (1+O(r^8))\int_{B(r)^k}\exp\Big\{-\int_{D(z)} K_0(u,u)\,du\Big\}\,dz. \end{aligned} \tag{65}$$

Moreover,

$$\begin{aligned} &\int_{B(r)^k}\exp\Big\{-\int_{D(z)} K_0(u,u)\,du\Big\}\,dz \\ &\quad = r^{2k}\int_{B(1)^k}\exp\Big\{-\int_{rD(z)} K_0(u,u)\,du\Big\}\,dz \\ &\quad = r^{2k}\int_{B(1)^k}\exp\Big\{-r^2\frac{1}{\pi}V(z)+\frac{r^2}{\pi}\int_{D(z)} e^{-|ru|^2}\,du\Big\}\,dz \\ &\quad = r^{2k}\int_{B(1)^k}\Big[1+\frac{r^2}{\pi}V(z)+O(r^4)\Big]\exp\Big\{-r^2\frac{1}{\pi}V(z)\Big\}\,dz. \end{aligned} \tag{66}$$

Formula (3) is a straightforward consequence of (62), (65), (66) and the asymptotic equality,

$$\frac{\int_{B(1)^k} V(z)\exp\{-r^2(1/\pi)V(z)\}\,dz}{\int_{B(1)^k}\exp\{-r^2(1/\pi)V(z)\}\,dz} = \frac{1}{\pi^k}\int_{B(1)^k} V(z)\,dz + O(r^2). \tag{67}$$

5.2. *Proof of Theorem 4, part* (b). For $z \in (\mathbb{C}\setminus B(R))^k$, let $\alpha_n$, $f_n$, $n\geq 1$, denote the eigenvalues and the eigenfunctions of the operator $K$ on $D(z)$ where $\alpha_1$ is the largest eigenvalue.

By Lemma 7,

$$\frac{R_{D(z)}(0,0)}{K(0,0)} = \frac{\sum_{n\geq 1}(\alpha_n/(1-\alpha_n))|f_n(0)|^2}{\sum_{n\geq 1}\alpha_n|f_n(0)|^2} \leq \frac{1}{1-\alpha_1}. \tag{68}$$

One has

$$\sum_{n\geq 1}\alpha_n=(1/\pi)V(z). \tag{69}$$

Note also that when $z\in(\mathbb{C}\setminus B(R))^k$, there exists $a\in\mathbb{C}$, $|a|=R$ such that $D(z)\supset B(a,R)$. Thus

$$\begin{aligned}\sum_{n\geq 2}(\alpha_n)^2&=\int_{D(z)}K^{(2)}_{D(z)}(u,u)\,du=\frac{1}{\pi^2}\int_{D(z)^2}e^{-|u-v|^2}\,du\,dv\\ &\geq\frac{1}{\pi^2}\int_{B(a,R)^2}e^{-|u-v|^2}\,du\,dv=\frac{1}{\pi^2}\int_{B(R)^2}e^{-|u-v|^2}\,du\,dv.\end{aligned} \tag{70}$$

Introduce

$$(*)=\frac{R_{D(z)}(0,0)}{K(0,0)}P\{N_{D(z)}(\phi)=0\}.$$

Thus, by (14), (31), (68) and (69),

$$(*)\leq\exp\left\{\sum_{n\geq 2}\log(1-\alpha_n)\right\}\leq\exp\left[-\frac{V(z)}{\pi}+\alpha_1+\frac{(\alpha_1)^2}{2}-\frac{1}{2}\sum_{n\geq 1}(\alpha_n)^2\right],$$

and hence by (70), we obtain

$$(*)\leq\exp\left[-\frac{V(z)}{\pi}+\frac{3}{2}-\frac{1}{2\pi^2}\int_{B(R)^2}e^{-|u-v|^2}\,du\,dv\right]. \tag{71}$$

From (57), (71), notation of Theorem 4 and (60), it follows that

$$\begin{aligned}E[V^k(\mathcal{C}\setminus B(R))]&=\int_{(\mathbb{C}\setminus B(R))^k}\frac{R_{D(z)}(0,0)}{K(0,0)}P\{N_{D(z)}(\phi)=0\}\,dz\\ &\leq e^{(3/2)-J(R)}\int_{(\mathbb{C}\setminus B(R))^k}\exp\left[-\frac{V(z)}{\pi}\right]dz\\ &=EV^k(\mathcal{C}_p\setminus B(R))\cdot e^{(3/2)-J(R)},\end{aligned} \tag{72}$$

which proves part (b) of Theorem 4.

PROBLEM 4. The facts that $EV(\mathcal{C})=EV(\mathcal{C}_p)=\pi$ and that the Ginibre–Voronoi tessellation is more regular than the Poisson–Voronoi tessellation suggest the conjecture (which seems to be confirmed by Monte Carlo simulation [16]) that the inequality,

$$EV^2(\mathcal{C})\leq EV^2(\mathcal{C}_p),$$

holds. It would be interesting to provide a rigourous proof of this property.

5.3. *Proof of Theorem 5.* By Theorem 1 and formula (54), the typical cell $\mathcal{C}$ coincides in law with the cell $C(0)$ related to the process $\phi_0$ which is obtained by removing from $\phi$ the point $Z$. If we consider the cell,

$$C_0 = \{z \in \mathbb{C}; \forall v \in \phi, |z| \leq |z - v|\} \subset \mathbb{C},$$

that is, if we do not remove the point $Z$, then for every Borel set $A \subset \mathbb{C}$ and $k \geq 1$,

$$\begin{aligned} E[V^k(C_0 \cap A)] &= \int_{A^k} P\{N_{D(z)}(\phi) = 0\}\, dz \\ &= \int_{A^k} \exp\left\{-\frac{V(z)}{\pi} - \sum_{n\geq 2} \frac{1}{n} \int K_{D(z)}^{(n)}(u,u)\, du\right\} dz \qquad (73) \\ &\leq \int_{A^k} \exp\{-V(z)/\pi\}\, dz = E[V^k(\mathcal{C}_p \cap A)]. \end{aligned}$$

Say that a point $u \in \phi$ is a neighbor of the origin if the bisecting line of the segment $[0, u]$ intersects the boundary of the cell $C_0$.

Denote by $\mathcal{N} = \mathcal{N}(\phi)$ the set of neighbours of the origin. Recall property (56), that is,

$$z \in C(0)^k \iff D(z) \cap \phi_0 = \varnothing.$$

By Theorem 1, we have also

$$\begin{aligned} z \in C_0^k &\iff D(z) \cap \phi = \varnothing \\ &\iff D(z) \cap \phi_0 = \varnothing \quad \text{and} \quad Z \notin D(z). \end{aligned} \qquad (74)$$

Moreover, if $Z \notin \mathcal{N}$ then obviously $C(0) = C_0$. Consequently, we obtain

$$\begin{aligned} &E[V^k(C(0)) - V^k(C_0)] \\ &\quad = E[\{V^k(C(0)) - V^k(C_0)\} \times 1_{\{Z\in\mathcal{N}\}}] \\ &\quad = \int_{\mathbb{C}^k} [P\{D(z) \cap \phi_0 = \varnothing, Z \in \mathcal{N}\} \\ &\qquad\quad - P\{D(z) \cap \phi_0 = \varnothing, Z \notin D(z), Z \in \mathcal{N}\}]\, dz \qquad (75) \\ &\quad = \int_{\mathbb{C}^k} P\{D(z) \cap \phi_0 = \varnothing, Z \in D(z), Z \in \mathcal{N}\}\, dz \\ &\quad \leq P\{Z \in \mathcal{N}\}^{1/2} \int_{\mathbb{C}^k} P\{D(z) \cap \phi_0 = \varnothing, Z \in D(z)\}^{1/2}\, dz. \end{aligned}$$

Now by (53), (19), (31) and notation (5),

$$P\{D(z) \cap \phi_0 = \varnothing, Z \in D(z)\}$$

$$
\begin{aligned}
&= \left(1 - \frac{K(0,0)}{R_{D(z)}(0,0)}\right) P\{N_{D(z)}(\phi_0) = 0\} \\
&= (\pi R_{D(z)}(0,0) - 1) P\{N_{D(z)}(\phi) = 0\} \\
&= H(z).
\end{aligned} \tag{76}
$$

Moreover, by (19),

$$
\begin{aligned}
&E[V^k(C(0)) - V^k(C_0)] \\
&\quad = \int_{\mathbb{C}^k} (P\{N_{D(z)}(\phi_0) = 0\} - P\{N_{D(z)}(\phi) = 0\})\, dz \\
&\quad = \int_{\mathbb{C}^k} H(z)\, dz.
\end{aligned} \tag{77}
$$

Relations (75)–(77) imply (6).

5.4. *Theorem 5, case $k = 1$.* Fix $a \in \mathbb{C}$ and consider the kernel $\widehat{K}(z_1, z_2) = K(z_1 - a, z_2 - a)$; hence

$$
\begin{aligned}
\widehat{K}(z_1, z_2) = \sum_{n \ge 0} & \frac{e^{-(1/2)|z_1 - a|^2} (z_1 - a)^n}{\sqrt{\pi n!}} \\
& \times \frac{e^{-(1/2)|z_2 - a|^2} (\overline{z_2 - a})^n}{\sqrt{\pi n!}}.
\end{aligned} \tag{78}
$$

Observe that the functions

$$
f_n(z) = e^{-(1/2)|z-a|^2} (z - a)^n, \qquad z \in D(a), n \ge 1,
$$

are orthogonal on $D(a) = B(a, r)$ with $r = |a|$ and that

$$
\int_{D(a)} |f_n(z)|^2\, dz = \pi \gamma(n + 1, r^2),
$$

where

$$
\gamma(n, u) = \int_0^u e^{-t} t^{n-1}\, dt
$$

is the incomplete gamma function.

Denote

$$
\alpha_n = \frac{\gamma(n + 1, r^2)}{n!}, \qquad n \ge 0, \tag{79}
$$

and

$$
\widehat{f}_n = (\pi \gamma(n + 1, r^2))^{-1/2} f_n,
$$

then

$$\widehat{K}(z_1, z_2) = \sum_{n\geq 0} \alpha_n \widehat{f}_n(z_1)\overline{\widehat{f}_n(z_2)}. \tag{80}$$

It follows from (80) that $\alpha_n$, $n \geq 0$, are the eigenvalues of the integral kernel $\widehat{K}$ on $D(a) = B(a, r)$ with $r = |a|$ and that for the resolvent kernel $\widehat{R}_{D(a)}$, we have

$$\begin{aligned} \widehat{R}_{D(a)}(0,0) &= \sum_{n\geq 0} \frac{\alpha_n}{1-\alpha_n} \widehat{f}_n(0)\overline{\widehat{f}_n(0)} \\ &= \frac{1}{\pi} \sum_{n\geq 0} \frac{r^{2n}e^{-r^2}}{\Gamma(n+1, r^2)}, \end{aligned} \tag{81}$$

where

$$\Gamma(n, u) = \Gamma(n) - \gamma(n, u) = \int_u^{+\infty} e^{-t}t^{n-1}\, dt.$$

By Remark 17 and formulas (14), (31) and (79),

$$P\{N_{D(a)}(\phi) = 0\} = \prod_{n\geq 0} \frac{\Gamma(n+1, r^2)}{n!}, \qquad r = |a|. \tag{82}$$

Therefore, by (81) and (82),

$$\begin{aligned} EV(\mathcal{C}) &= \int_{\mathbb{C}} P\{N_{D(z)}(\phi_0) = 0\}\, dz \\ &= \int_{\mathbb{C}} P\{N_{D(z)}(\phi) = 0\} \frac{\widehat{R}_{D(z)}(0,0)}{\widehat{K}_{D(z)}(0,0)}\, dz \\ &= \pi \int_0^{+\infty} \prod_{n\geq 0} \frac{\Gamma(n+1,t)}{n!} \times \sum_{n\geq 0} \frac{t^n e^{-t}}{\Gamma(n+1,t)}\, dt \\ &= \pi \left[ - \prod_{n\geq 0} \frac{\Gamma(n+1,t)}{n!} \right]_{t=0}^{t=+\infty} = \pi. \end{aligned}$$

This is the expected result. Now, we have also

$$\begin{aligned} EV(C_0) &= \int_{\mathbb{C}} P\{N_{D(z)}(\phi) = 0\}\, dz \\ &= \pi \int_0^{+\infty} \prod_{n\geq 0} \frac{\Gamma(n+1,t)}{n!}\, dt. \end{aligned}$$

In [20] [formula (15.1.27)], M. L. Mehta showed that, for every $t \geq 0$,

$$\prod_{n\geq 0} \frac{\Gamma(n+1,t)}{n!} \leq (1+t)e^{-2t}. \tag{83}$$

This implies

$$EV(C_0) \leq \frac{3\pi}{4}.$$

From (79), (81), and with the notation (5), we obtain

$$\begin{aligned} H(a) &= (\pi \widehat{R}_{D(a)}(0,0) - 1) \prod_{n\geq 0}(1-\alpha_n) \\ &= \prod_{n\geq 0} \frac{\Gamma(n+1,r^2)}{n!} \times \sum_{n\geq 0} \frac{r^{2n}e^{-r^2}}{\Gamma(n+1,r^2)} \\ &\quad - \prod_{n\geq 0} \frac{\Gamma(n+1,r^2)}{n!} \\ &= \prod_{n\geq 0} \frac{\Gamma(n+1,r^2)}{n!} \times \sum_{n\geq 1} \frac{r^{2n}e^{-r^2}}{\Gamma(n+1,r^2)}. \end{aligned} \tag{84}$$

Consequently,

$$\int_{\mathbb{C}} H(z)\,dz = 2\pi\left(1 - \int_0^{+\infty} \prod_{n\geq 0} \frac{\Gamma(n+1,t)}{n!}\,dt\right) \tag{85}$$

and inserting (84) and (85) in (6), we obtain (7). Now, applying the Hölder inequality, we get

$$\begin{aligned} &\left(\int_0^{+\infty} \sqrt{\left(\sum_{n\geq 1} \frac{t^n e^{-t}}{\Gamma(n+1,t)}\right) \prod_{n\geq 0} \frac{\Gamma(n+1,t)}{n!}}\,dt\right)^2 \\ &\quad \leq \int_0^{+\infty} \left(\sum_{n\geq 1} \frac{t^n e^{-t}}{\Gamma(n+1,t)}\right) \prod_{n\geq 1} \frac{\Gamma(n+1,t)}{n!}\,dt \\ &\quad \times \int_0^{+\infty} e^{-t}\,dt = 1. \end{aligned}$$

This and inequality (7) give

$$P\{Z \in \mathcal{N}(\phi)\} \geq \left[1 - \int_0^{+\infty} \prod_{n\geq 0} \frac{\Gamma(n+1,t)}{n!}\,dt\right]^2$$

and by (83) we obtain

$$P\{Z \in \mathcal{N}(\phi)\} \geq \tfrac{1}{16}.$$

PROBLEM 5. It would be interesting to investigate other geometric characteristics of $\mathcal{C}$, among others, the number of sides, the perimeter, and the radius of the smaller disc containing $C(0)$.

## APPENDIX: PROOF OF LEMMA 13

Let $\alpha$, $\beta$ be point processes on $\mathbb{E}$ such that for every decreasing event $\mathcal{A} \in \mathcal{F}_d(\mathbb{E})$,

$$P\{\alpha \in \mathcal{A}\} \leq P\{\beta \in \mathcal{A}\}. \tag{86}$$

Recall that $\mathcal{F}_d \subset \mathcal{F}$ is the collection of sets which are a finite union of elementary decreasing events. We want to prove that point process $\alpha \in \mathcal{M}_\sigma(\mathbb{E})$ stochastically dominates the point process $\beta \in \mathcal{M}_\sigma(\mathbb{E})$ which is equivalent to the fact that inequality (86) above is satisfied for every decreasing event $\mathcal{A} \in \mathcal{F}(\mathbb{E})$.

**A.1. Step I.** It suffices to prove that (86) is satisfied for every decreasing event $\mathcal{A} \in \mathcal{F}(\mathbb{E})$ by assuming that $\mathcal{A} \in \mathcal{F}(\mathbb{E})$ is a closed subset of $\mathcal{M}_\sigma(\mathbb{E})$. Moreover, if $\mathbb{E}$ is compact then we may assume that $\mathcal{A} \in \mathcal{F}(\mathbb{E})$ is compact as well.

Indeed, denote by $Q$ the law of point process $\alpha$ and by Q' the law of the process $\beta$. We want to show that $Q(\mathcal{A}) \leq Q'(\mathcal{A})$ for all $\mathcal{A} \in \mathcal{F}(\mathbb{E})$. Recall that $\mathcal{M}_\sigma(\mathbb{E})$ is a Polish space. Then, by the Lusin theorem, $Q(\mathcal{A}) = \sup\{Q(\mathfrak{A}); \mathfrak{A} \subseteq \mathcal{A}$ and $\mathfrak{A}$ is compact$\}$ and $Q'(\mathcal{A}) = \sup\{Q'(\mathfrak{A}); \mathfrak{A} \subseteq \mathcal{A}$ and $\mathfrak{A}$ is compact$\}$. Consider now a compact set $\mathfrak{A} \subseteq \mathcal{A}$ and denote

$$\tilde{\mathfrak{A}} = \{\xi \in \mathcal{M}_\sigma(\mathbb{E}); \exists \zeta \in \mathfrak{A} \text{ such that } \xi \subseteq \zeta\}.$$

The set $\tilde{\mathfrak{A}}$ is decreasing and $\mathfrak{A} \subseteq \tilde{\mathfrak{A}} \subseteq \mathcal{A}$. Consequently the result follows from the lemma below.

LEMMA 26. (i) *The set* $\tilde{\mathfrak{A}}$ *above is closed.*
(ii) *If* $\mathbb{E}$ *is compact then the set* $\tilde{\mathfrak{A}}$ *is compact as well.*

PROOF. For property (i), consider a sequence $(\xi_n)_n \subset \tilde{\mathfrak{A}}$ such that $\xi_n \to \xi \in \mathcal{M}_\sigma(\mathbb{E})$ [the space $\mathcal{M}_\sigma(\mathbb{E})$ being endowed with the vague topology]. We want to show that $\xi \in \tilde{\mathfrak{A}}$. For every $n \geq 1$ there exists $\zeta_n \in \mathfrak{A}$ such that $\xi_n \subseteq \zeta_n$. With the set $\mathfrak{A}$ being compact, there exists a convergent subsequence $\zeta_{n_k} \to \zeta \in \mathfrak{A}$. We claim that $\xi \subseteq \zeta$ (and thus $\xi \in \tilde{\mathfrak{A}}$). Indeed, suppose that there exists $x \in \xi$ such that $x \notin \zeta$. Let $f_i : \mathbb{E} \to \mathbb{R}$, $i = 1, 2$, be continuous functions with compact supports, respectively, $K_i$, $i = 1, 2$, such that $K_1 \subset K_2$, $\xi \cap K_1 = x$, $K_2 \cap \zeta = \varnothing$, $0 \le f_1, f_2 \le 1$, $f_1(x) = 1$ and $f_2 \equiv 1$ on $K_1$. Denote $V(\xi) = \{\eta; |1 - \sum_{z\in\eta} f_1(z)| \le 1/2\}$ and $V(\zeta) = \{\eta; |\sum_{z\in\eta} f_2(z)| \le 1/2\}$. For large $n_k$ we have $\xi_{n_k} \in V(\xi)$ and $\zeta_{n_k} \in V(\zeta)$ from which follows that $\xi_{n_k} \cap K_1 \neq \varnothing$ and $\zeta_{n_k} \cap K_1 = \varnothing$ which implies in turn the contradiction $\xi_{n_k} \nsubseteq \zeta_{n_k}$. Property (i) is then proved. To prove property (ii) notice (see [14]) that the compactness of the sets $\mathbb{E}$ and $\mathfrak{A}$ implies that there exists $A > 0$ such that $\mathfrak{A} \subset \{N_{\mathbb{E}} \le A\}$. Obviously $\tilde{\mathfrak{A}}$ is included in $\{N_{\mathbb{E}} \le A\}$, the later being compact (see [14], page 33); the result follows from property (i). □

**A.2. Step II.** We may suppose that $\mathbb{E}$ is a compact separable metric space. Indeed, we have $\mathbb{E} = \bigcup_{n\ge1} K_n$ where the sets $K_n \uparrow \mathbb{E}$ are compact with countable bases. Denote $\alpha_n = \alpha \cap K_n$ and $\beta_n = \beta \cap K_n$. Condition (86) implies that $P\{\alpha_n \in \mathcal{A}\} \le P\{\beta_n \in \mathcal{A}\}$ is satisfied for every decreasing event $\mathcal{A} \in \mathcal{F}_d(K_n)$. Thus if condition (86) implies stochastic domination for $\mathbb{E}$ compact, then the process $\beta_n$ is stochastically dominated by the process $\alpha_n$, and we have $P\{\alpha_n \in \mathcal{A}\} \le P\{\beta_n \in \mathcal{A}\}$ for $\mathcal{A} \in \mathcal{F}(\mathbb{E})$. This and the fact that $\alpha_n \uparrow \alpha$ and $\beta_n \uparrow \beta$ implies that $P\{\alpha \in \mathcal{A}\} \le P\{\beta \in \mathcal{A}\}$ for closed decreasing sets $\mathcal{A} \in \mathcal{F}(\mathbb{E})$.

**A.3. Step III.** We suppose that $\mathbb{E}$ is compact with a metric $d$. Fix a compact decreasing set $\mathcal{A} \in \mathcal{F}(\mathbb{E})$. We can suppose that there exists $\varepsilon > 0$ such that for each $A \in \mathcal{A}$ and each $x, y \in A$, we have $d(x, y) > \varepsilon$. Indeed, denote by $\mathcal{A}_n \subset \mathcal{A}$ the set where the elements are the finite sets $A \in \mathcal{A}$ such that for every $x, y \in A$, we have $d(x, y) > 1/n$. Then the sets $\mathcal{A}_n$ are decreasing as well and when $n \to +\infty$, we have $Q(\mathcal{A}_n) \to Q(\mathcal{A})$ and $Q'(\mathcal{A}_n) \to Q'(\mathcal{A})$.

**A.4. Step IV.** For each $A \in \mathcal{A}$ (note that $A$ is finite) and $n > 1/\varepsilon$ consider the set

$$O_{n,A} = \{\phi \in \mathcal{M}_\sigma(\mathbb{E}); N_{B_0(x,1/n)}(\phi) \le 1 \text{ for each } x \in A \\ \text{and } N_{E\setminus\bigcup_{x\in A} B_0(x,1/n)}(\phi) = 0\}, \tag{87}$$

where $B_0(x, 1/n)$ is the open ball.

Denote $\mathfrak{K}_n = \bigcup_{A\in\mathcal{A}} O_{n,A}$.

In order to finish it suffices to note that:

- The sets $O_{n,A}$ are open (see [14]);
- We have $\mathfrak{K}_{n+1} \subset \mathfrak{K}_n$;
- $\bigcap \mathfrak{K}_n = \mathcal{A}$;

– $\mathfrak{K}_n$ is a covering of $\mathcal{A}$ by open sets, $\mathcal{A}$ being compact there exists a covering of $\mathcal{A}$ by a finite number of sets $O_{n,A}$.

Consequently, in order to obtain $Q(\mathcal{A}) \leq Q'(\mathcal{A})$ it suffices to have $Q(\bigcup_{i=1,\dots,N} O_{n,A_i}) \leq Q'(\bigcup_{i=1,\dots,N} O_{n,A_i})$. Naturally, $\bigcup_{i=1,\dots,N} O_{n,A_i} \in \mathcal{F}_d$ which completes the proof.

Université Claude Bernard Lyon 1
Département de Mathématiques
Institut Camille Jordan UMR 5208 CNRS
bât. J. Braconnier
43, boulevard du 11 novembre 1918
69622 Villeurbanne Cedex
France
E-mail: Andre.Goldman@univ-lyon1.fr